\documentclass[12pt]{article}
\usepackage{epsfig}
\usepackage{amssymb}
\usepackage{graphicx}
\usepackage{amsmath}
\usepackage{amsthm}
\usepackage{amscd}
\usepackage{arydshln}
\usepackage{mathtools}
\usepackage{units}
\usepackage[all,cmtip]{xy}
\usepackage{diagxy}
\usepackage{xymatrix}
\usepackage[all]{xy}
\usepackage{fancyhdr}
\usepackage[nottoc,numbib]{tocbibind}
\usepackage{graphicx}
\usepackage{mathrsfs}
\usepackage[perpage]{footmisc}
\usepackage{titlesec}
\usepackage{color}
\usepackage{geometry}
\usepackage{multicol}
\usepackage{bbold}
\usepackage{slashed}
\usepackage{changepage}
\usepackage[square,numbers]{natbib}
\usepackage[bookmarks=false]{hyperref}

\newcommand{\tens}[1]{\mathbin{\mathop{\otimes}\limits_{#1}}}

\newtheoremstyle{basic}{11pt}{11pt}{}{}{\bfseries}{.}{0.5em}{}
\newtheoremstyle{proof}{11pt}{11pt}{}{}{\scshape}{:}{0.5em}{}

\newtheorem{prop}{Proposition}
\newtheorem{cor}{Corollary}

\newtheorem{thm}{Theorem}

\theoremstyle{basic}
\newtheorem{defn}{Definition}

\theoremstyle{basic}
\newtheorem{exa}{Example}

\theoremstyle{basic}
\newtheorem{rem}{Remark}

\theoremstyle{proof}
\newtheorem*{prf}{Proof}

\makeatletter
\newcommand{\colim@}[2]{%
  \vtop{\m@th\ialign{##\cr
    \hfil$#1\operator@font colim$\hfil\cr
    \noalign{\nointerlineskip\kern1.5\ex@}#2\cr
    \noalign{\nointerlineskip\kern-\ex@}\cr}}%
}
\newcommand{\colim}{%
  \mathop{\mathpalette\colim@{\rightarrowfill@\textstyle}}\nmlimits@
}
\makeatother

\titleformat*{\section}{\large\bfseries}
\titleformat*{\subsection}{\normalsize\bfseries}

\begin{document}

\title{Hopf algebroids, Hopf categories,\\
and their Galois theories}

\author{Clarisson Rizzie Canlubo\\
University of Copenhagen\\
\textit{clarisson@math.ku.dk}}

\maketitle



\pagestyle{fancy}
\headheight=17pt
\textwidth=500pt
\fancyhead{}
\fancyhead[L]{\footnotesize{\leftmark}}
\fancyhead[R]{\thepage}
\renewcommand{\headrulewidth}{1pt}
\renewcommand{\footrulewidth}{1pt}
\fancyfoot{}
\fancyfoot[R]{}
\fancyfoot[L]{\scriptsize{Hopf algebroids, Hopf categories, and their Galois theories}}
\fancyfoot[C]{}


\begin{abstract}
Hopf algebroids are generalization of Hopf algebras over non-commutative base rings. It consists of a left- and a right-bialgebroid structure related by a map called the antipode. However, if the base ring of a Hopf algebroid is commutative one does not necessarily have a Hopf algebra. Meanwhile, a Hopf category is the categorification of a Hopf algebra. It consists of a category enriched over a braided monoidal category such that every hom-set carries a coalgebra structure together with an antipode functor. In this article, we will introduce the notion of a topological Hopf category$-$ a small category whose set of objects carries a topology and whose categorical structure maps are sufficiently continuous. The main result of this paper is to describe the relation between finitely-generated projective Hopf algebroids over commutative unital $C^{*}$-algebras and topological coupled Hopf categories of finite-type whose space of objects is compact and Hausdorff. We will accomplish this by using methods in algebraic geometry and spectral theory. Lastly, we will show that not only the two objects are tightly related, but so are their respective Galois theories.

\


\

\end{abstract}

\tableofcontents

\normalsize

\section{Introduction}\label{S1.0}

Hopf algebras are robust generalization of groups. Recently, many authors have studied much more general Hopf-like structures: weak Hopf algebras, Hopf monads, $\times_{R}$-Hopf algebras, compact quantum groups to name a few. In this article, we will mainly be interested with Hopf algebroids and Hopf categories. In the literature, there are plenty of inequivalent notions of a Hopf algebroid. For the exposition on these notions, see B\"ohm \cite{bohm}. Batista et al. \cite{bcv} introduced the notion of a Hopf category which is the natural categorification of a Hopf algebra. Motivated by a fundamental related to his PhD thesis, the author tries to describe the geometry of Hopf algebroids over $C(X)$. This geometric description necessitates a structure closely related to a Hopf category, but which has not appeared in the literature as far as the author's knowledge. We will define such structures in section (\ref{S3.0}).

We will recall in section (\ref{S2.0}) definitions and properties of Hopf algebroids. For completeness, we will also give a short exposition on the representation theoretic and Galois theoretic aspects of Hopf algebroids. Most of section (\ref{S2.0}) follows \cite{bohm} except for the definition of morphisms of Hopf algebroids and the definition of a coupled Hopf algebra. In section (\ref{S3.0}), we will define what topological Hopf categories are and we will also define coupled Hopf categories and their topological version. We will end that section with a formulation of Galois theory for Hopf categories and all the related variant we will introduce in that section.

One of the main result of this paper is theorem (\ref{T4.1}). It gives a bijective correspondence between finitely-generated projective Hopf algebroids over $C(X)$ and topological coupled Hopf categories of finite type. Using algebraic geometric and spectral theoretic methods, spanning the entirety of section (\ref{S4.0}), we will prove this result. The second main result is theorem (\ref{T5.1}), which states that, not only is there a bijection between Hopf algebroids and topological Hopf categories, their Galois theories also matched in a bijective manner.

Following David Hilbert's statement:

\vspace{.1in}

\begin{adjustwidth}{2cm}{2cm}
\textit{"The art of doing mathematics consists in finding that special case which contains all the germs of generality."}
\end{adjustwidth}

\vspace{.1in}

\noindent we will discuss a very important example in section (\ref{S3.2}) which completely illustrates the general situation.


\textbf{Acknowledgement.} I would like to thank my PhD supervisor Ryszard Nest for guiding me through my studies in non-commutative geometry and for the valuable discussions that help me write this article. I would also like to thank DSF Grant, UP Diliman and the support of the Danish National Research Foundation through the Centre for Symmetry and Deformation (DNRF92).


\section{Hopf algebroids}\label{S2.0}

\subsection{Definitions}\label{S2.1}

There are several inequivalent notions of a Hopf algebroid. We will briefly present here the one defined in B\"ohm \cite{bohm}. An $R$-\textit{ring} is a monoid object in the category of $R$-bimodules. Explicitly, an $R$-ring is a triple $(A,\mu,\eta)$ where $A\otimes_{R}A\stackrel{\mu}{\longrightarrow}A$ and $R\stackrel{\eta}{\longrightarrow}A$ are $R$-bimodule maps satisfying the associativity and unit axioms similar for algebras over commutative rings. A morphism of $R$-rings is a monoid morphism in category of $R$-bimodules. It is important to note that there is a bijection between $R$-rings $(A,\mu,\eta)$ and $k$-algebra morphisms $R\stackrel{\eta}{\longrightarrow}A$. Similar to the case of algebras over commutative rings, we can define modules over $R$-rings. For an $R$-ring $(A,\mu,\eta)$, a \textit{right} (resp. \textit{left}) $(A,\mu,\eta)$-\textit{module} is an algebra for the monad $-\otimes_{R}A$ (resp. $A\otimes_{R}-$) on the category $\mathcal{M}_{R}$ (resp. ${}_{R}\mathcal{M}$) of right (resp. left) modules over $R$.

We can dualize all the objects we have defined in the previous paragraph. An $R$-\textit{coring} is a comonoid in the category of $R$-bimodules, i.e a triple $(C,\Delta,\epsilon)$ where $C\stackrel{\Delta}{\longrightarrow}C\otimes_{R}C$ and $C\stackrel{\epsilon}{\longrightarrow}R$ are $R$-bimodule maps satisfying the coassociativity and counit axioms dual to those axioms satisfied by the structure maps of an $R$-ring. A morphism of $R$-corings is a morphism of comonoids. Given an $R$-coring $(C,\Delta,\epsilon)$, similar to coalgebras over commutative rings, we define a \textit{right} (resp. \textit{left}) $(C,\Delta,\epsilon)$-\textit{comodule} as a coalgebra for the comonad $-\otimes_{R}C$ (resp. $C\otimes_{R}-$) on the category $\mathcal{M}_{R}$ (resp. ${}_{R}\mathcal{M}$).

\begin{defn}\label{D2.1}
A \textit{right} (resp. \textit{left}) $R$-\textit{bialgebroid} $B$ is an $R\otimes_{k}R^{op}$-ring $(B,s,t)$ and an $R$-coring $(B,\Delta,\epsilon)$ satisfying:

\begin{enumerate}
\item [(a)] $R\stackrel{s}{\longrightarrow}B$ and $R^{op}\stackrel{t}{\longrightarrow}B$ are $k$-algebra maps with commuting images defining the $R\otimes_{k}R^{op}$-ring structure on $B$ which is compatible to the $R$-bimodule structure as an $R$-coring thru the following relation:
\[ r\cdot b \cdot r':=bs(r')t(r), \ \ (\text{resp.} \ r\cdot b \cdot r':=s(r)t(r')b,) \hspace{.15in} \forall r,r'\in R, b\in B. \]

\item [(b)] With the above $R$-bimodule structure on $B$ one can form $B\otimes_{R}B$. The coproduct $\Delta$ is required to corestrict to a $k$-algebra map to
\[ B\times_{R}B:=\left\{\sum\limits_{i}b_{i}\otimes_{R}b_{i}'\left|\sum\limits_{i}s(r)b_{i}\otimes_{R}b_{i}'=\sum\limits_{i}b_{i}\otimes_{R}t(r)b_{i}',\forall r\in R\right.\right\} \]
\noindent respectively,
\[ B \prescript{}{R}{\times} \ B:=\left\{\sum\limits_{i}b_{i}\otimes_{R}b_{i}'\left|\sum\limits_{i}b_{i}t(r)\otimes_{R}b_{i}'=\sum\limits_{i}b_{i}\otimes_{R}b_{i}'s(r),\forall r\in R\right.\right\}. \]

\item [(c)] The counit $B\stackrel{\epsilon}{\longrightarrow}R$ extends the right (resp. left) regular $R$-module structure on $R$ to a right (resp. left) $(B,s)$-module.
\end{enumerate}

\noindent A \textit{morphism} of $R$-bialgebroids is a morphism of $R\otimes R^{op}$-rings and $R$-corings.
\end{defn}

\begin{rem}\label{R2.1}
\begin{enumerate}
\item[]
\item[(1)] The $k$-algebra maps $s$ and $t$ define a $k$-algebra map $\eta=s\otimes_{k}t$. As we have noted, such $k$-algebra uniquely determines an $R\otimes_{k}R^{op}$-ring structure on $B$. The maps $s$ and $t$ are called the \textit{source} and \textit{target} maps, respectively.

\item[(2)] The $k$-submodule $B\times_{R}B$ (resp. $B\prescript{}{R}{\times} \ B$) of $B\otimes_{R}B$ is a $k$-algebra with factorwise multiplication. This is called the \textit{Takeuchi product}. The map $R\otimes_{k}R^{op}\longrightarrow B\times_{R}B$, $r\otimes_{k}r'\mapsto t(r')\otimes_{R}s(r)$ is easily seen to be a $k$-algebra morphism and hence, $B\times_{R}B$ is an $R\otimes_{k}R^{op}$-ring. The corestriction of $\Delta$ is an $R\otimes_{k}R^{op}$-bimodule map. Hence, $\Delta$ is an $R\otimes R^{op}$-ring map. The same is true for $B\prescript{}{R}{\times} \ B$.

\item[(3)] The source map $s$ is a k-algebra map and so it defines a unique $R$-ring structure on $B$. The right version of condition (c) explicitly means that $r\cdot b:=\epsilon(s(r)b)$, $\forall r\in R, b\in B $ defines a right $(B,s)$-action on $R$.
\end{enumerate}
\end{rem}

\begin{defn}\label{D2.2}
Let $k$ be a commutative, associative unital ring and let $L$ and $R$ be associative $k$-algebras. A \textit{Hopf algebroid} $\mathcal{H}$ is a triple $\mathcal{H}=(H_{L},H_{R},S)$. $H_{L}$ and $H_{R}$ are bialgebroids having the same underlying $k$-algebra $H$. Specifically, $H_{L}$ is a left $L$-bialgebroid with $(H,s_{L},t_{L})$ and $(H,\Delta_{L},\epsilon_{L})$ as its underlying $L\otimes_{k}L^{op}$-ring and $L$-coring structures. Similarly, $H_{R}$ is a right $R$-bialgebroid with $(H,s_{R},t_{R})$ and $(H,\Delta_{R},\epsilon_{R})$ as its underlying $R\otimes_{k}R^{op}$-ring and $R$-coring structures. Let us denote by $\mu_{L}$ (resp. $\mu_{R}$) the multiplication on $(H,s_{L})$ (resp. $(H,s_{R})$). $S$ is a (bijective) $k$-module map $H\stackrel{S}{\longrightarrow}H$, called the \textit{antipode}. The compatibility conditions of these structures are as follows.

\begin{enumerate}
\item[(a)] the sources $s_{R},s_{L}$, targets $t_{R},t_{L}$ and counits $\epsilon_{R},\epsilon_{L}$ satisfy

\[ s_{L}\circ \epsilon_{L}\circ t_{R}=t_{R}, \hspace{.15in} t_{L}\circ \epsilon_{L}\circ s_{R}=s_{R}, \hspace{.15in} s_{R}\circ \epsilon_{R}\circ t_{L}=t_{L}, \hspace{.15in} t_{R}\circ \epsilon_{R}\circ s_{L}=s_{L}, \]

\item[(b)] the left- and right-regular comodule structures commute, i.e.

\[ \xymatrix{
H \ar[d(1.7)]_-{\Delta_{L}} \ar[r(1.7)]^-{\Delta_{R}}
& & H\tens{R} H \ar[d(1.6)]^-{\Delta_{L}\tens{R} id}\\
& & \\
H\tens{L} H \ar[r(1.6)]_-{id \tens{L}\Delta_{R}}
& & H\tens{L} H\tens{R} H } \hspace{0.75in}
\xymatrix{
H \ar[d(1.7)]_-{\Delta_{R}} \ar[r(1.7)]^-{\Delta_{L}}
& & H\tens{L} H \ar[d(1.6)]^-{\Delta_{R}\tens{L} id}\\
& & \\
H\tens{R} H \ar[r(1.6)]_-{id \tens{R}\Delta_{L}}
& & H\tens{R} H\tens{L} H } \]

\item[(c)] for all $l\in L, r\in R$ and for all $h\in H$ we have $S(t_{L}(l)ht_{R}(r))=s_{R}(r)S(h)s_{L}(l)$,

\item[(d)] $S$ is the convolution inverse of the identity map i.e., the following diagram commute

\[ \xymatrix{
& & H\tens{L} H \ar[rrrr]^-{S\tens{L} id} & & & & H\tens{L} H \ar[rrd]^-{\mu_{L}} & & \\
H \ar[rru]^-{\Delta_{L}}  \ar[rrrr]^-{\epsilon_{R}} & & & & R \ar[rrrr]^-{s_{R}} & & & & H \\
H \ar[rrd]_-{\Delta_{R}}  \ar[rrrr]_-{\epsilon_{L}} & & & & L \ar[rrrr]_-{s_{L}} & & & & H \\
& & H\tens{R} H \ar[rrrr]_-{id\tens{R} S} & & & & H\tens{R} H \ar[rru]_-{\mu_{R}} & & \\
} \]

\end{enumerate}
\end{defn}

\begin{rem}\label{R2.2}
\begin{enumerate}
\item[]
\item[(1)] In the constituent bialgebroids $H_{R}$ and $H_{L}$, the counits $\epsilon_{R}$ and $\epsilon_{L}$ extend the regular module structures on the base rings $R$ and $L$ to the $R$-ring $(H,s_{R})$ and to the $L$-ring $(H,s_{L})$, respectively. Equivalently, the counits extend the regular module structures on the base rings $R$ and $L$ to the $R^{op}$-ring $(H,t_{R})$ and to the $L^{op}$-ring $(H,t_{L})$. This particularly implies that the maps $s_{L}\circ \epsilon_{L}$, $t_{L}\circ \epsilon_{L}$, $s_{R}\circ \epsilon_{R}$ and $t_{R}\circ \epsilon_{R}$ are idempotents. This means that the images of $s_{R}$ and $t_{L}$ coincides in $H$. Same is true for the images of $s_{L}$ and $t_{R}$.

\item[(2)] Notice that for condition (b) to make sense, apart from being an $L$-bimodule map, $\Delta_{L}$ has to be an $R$-bimodule map. This is the case using remark (1). Similarly, $\Delta_{R}$ is an $L$-bimodule map.

\item[(3)] We can equip $H$ with two $(R,L)$-bimodule structures one using $t_{R}$ and $t_{L}$ and the other using $s_{R}$ and $s_{L}$. Condition (c) relates these two $(R,L)$-bimodules structures via the antipode $S$ which in turn makes the diagram in condition (d) defined.

\item[(4)] A most convenient way to summarize the property of the antipode of a Hopf algebra is to express it as the inverse of the identity map in the convolution algebra of endomorphisms of that Hopf algebra. For Hopf algebroids, the antipode is the inverse of the identity map in the appropriate category, called the \textit{convolution category} of $\mathcal{H}$. As before, $R$ and $L$ are $k$-algebras. Let $X$ and $Y$ be $k$-modules such that $X$ has an $R$-coring $(X, \Delta_{R}, \epsilon_{R})$ and an $L$-coring $(X,\Delta_{L},\epsilon_{L})$ structures and $Y$ has an $L\otimes_{k}R$-ring structure with multiplications $\mu_{R}:Y\otimes_{R}Y\longrightarrow Y$ and $\mu_{L}:Y\otimes_{L} Y\longrightarrow Y$. Define the convolution category $Conv(X,Y)$ to be the category with two objects labelled $R$ and $L$. For $I,J\in \left\{R,L\right\}$, a morphism $I\longrightarrow J$ is a $J$-$I$-bimodule map $X\longrightarrow Y$. For $I,J,K\in \left\{R,L\right\}$ and morphisms $J\stackrel{f}{\longrightarrow}I$ and $K\stackrel{g}{\longrightarrow}J$, we define the composition $f\ast g$ to be the following convolution

\[ f\ast g = \mu_{J} \circ (f\tens{J}g) \circ \Delta_{J}. \]

\noindent The antipode $S$ of a Hopf algebroid $\mathcal{H}$ is the inverse of the identity map $H\stackrel{id}{\longrightarrow}H$ viewed as an arrow in $Conv(H,H)$.

\item[(5)] Let us note that condition (c) in the definition of a bialgebroid implies that $\epsilon_{L}\circ s_{L}:L\longrightarrow L$ is the identity. Similarly, $\epsilon_{R}\circ s_{R}:R\longrightarrow R$ is also the identity. Using condition (a) in the definition of a Hopf algebroid, we see that the following compositions define pairs of inverse $k$-algebra maps.

\[ \xymatrix{
L \ar[rr]^-{\epsilon_{R}\circ s_{L}} & & R^{op} \ar[rr]^-{\epsilon_{L}\circ t_{R}} & & L} \hspace{.5in}
\xymatrix{
R \ar[rr]^-{\epsilon_{L}\circ s_{R}} & & L^{op} \ar[rr]^-{\epsilon_{R}\circ t_{L}} & & R}\]

\noindent This is particular implies that $R$ and $L$ are anti-isomorphic $k$-algebras.

\item[(6)] Since there are two coproducts involved in a Hopf algebroid, namely $\Delta_{L}$ and $\Delta_{R}$, we will use different Sweedler notations for their corresponding components. We will write $\Delta_{L}(h)=h_{[1]}\otimes_{L}h_{[2]}$ and $\Delta_{R}(h)=h^{[1]}\otimes_{R}h^{[2]}$ for $h\in H$.

\item[(7)] With a fixed bijective antipode $S$, the constituent left- and right-bialgebroids of a Hopf algebroid determine each other, see for example the article \cite{bohm-szlachanyi}. In view of this and the fact that $L$ and $R$ are anti-isomorphic, in the sequel where we will be mainly interested with Hopf algebroids with bijective antipodes we will simply call $\mathcal{H}$ a Hopf algebroid \textit{over} $R$ instead of explicitly mentioning $L$.

\end{enumerate}
\end{rem}

\begin{defn} \label{D2.3}
Let $(H_{L},H_{R},S)$ and $(H_{L}^{'},H_{R}^{'},S^{'})$ be Hopf algebroids over $R$. An \textit{algebraic morphism} $(H_{L},H_{R},S)\longrightarrow(H_{L}^{'},H_{R}^{'},S^{'})$ of Hopf algebroids is a pair $(\varphi_{L},\varphi_{R})$ of a left-bialgebroid morphism $\varphi_{L}$ and a right-bialgebroid morphism $\varphi_{R}$ for which the following diagrams commute

\[ \xymatrix{
H_{L} \ar[dd]_-{\varphi_{L}} \ar[rr]^-{S}
& & H_{R} \ar[dd]^-{\varphi_{R}}\\
& & \\
H_{L}^{'} \ar[rr]_-{S^{'}}
& & H_{R}^{'} } \hspace{0.75in}
\xymatrix{
H_{R} \ar[dd]_-{\varphi_{R}} \ar[rr]^-{S}
& & H_{L} \ar[dd]^-{\varphi_{L}}\\
& & \\
H_{R}^{'} \ar[rr]_-{S^{'}}
& & H_{L}^{'} } \]

\noindent and composition of such a pair is componentwise.

Let $R$ and $R^{'}$ be $k$-algebras and $(H_{L},H_{R},S)$ and $(K_{L^{'}},K_{R^{'}},S^{'})$ be Hopf algebroids over $R$ and $R^{'}$, respectively. In view of remark (\ref{R2.2}) (7) above, denote by $L=R^{op}$ and $L^{'}=(R^{'})^{op}$. A \textit{geometric morphism} $(H_{L},H_{R},S)\longrightarrow(K_{L^{'}},K_{R^{'}},S^{'})$ of Hopf algebroids is a pair $(f,\phi)$ of $k$-algerba maps $R\stackrel{f}{\longrightarrow}R^{'}$ and $H\stackrel{\phi}{\longrightarrow} K$, where $H,K$ denote the underlying $k$-algebra structures of the Hopf algebroids under consideration. These two maps satisfy the following compatibility conditions.

\begin{enumerate}

\item[(a)] $f$ and $\phi$ intertwines the source, target and counit maps of the left-bialgebroid structures of $\mathcal{H}$ and $\mathcal{K}$, i.e.

\[ \xymatrix{
H \ar[dd]_-{\phi} \ar[rr]^-{\epsilon^{H}_{L}}
& & L \ar[dd]^-{f}\\
& & \\
K \ar[rr]_-{\epsilon^{K}_{L}}
& & L^{'} } \hspace{0.15in}
\xymatrix{
L \ar[dd]_-{f} \ar[rr]^-{t^{H}_{L}}
& & H \ar[dd]^-{\phi}\\
& & \\
L^{'} \ar[rr]_-{t^{K}_{L}}
& & K } \hspace{0.15in}
\xymatrix{
L \ar[dd]_-{f} \ar[rr]^-{s^{H}_{L}}
& & H \ar[dd]^-{\phi}\\
& & \\
L^{'} \ar[rr]_-{s^{K}_{L}}
& & K. } \]

\noindent Same goes for the source, target and counit maps of the right-bialgebroid structures.

\item[(b)] In view of condition $(a)$, the $k$-bimodule map $\phi\otimes_{k}\phi$ defines $k$-bimodule maps

\[  \xymatrix{ H \prescript{}{L}{\otimes} \ H \ar[rr]^-{\phi\prescript{}{f}{\otimes} \ \phi} && K \prescript{}{L^{'}}{\otimes} \ K, } \hspace{.25in} \xymatrix{ H \otimes_{R} \ H \ar[rr]^-{\phi\ \otimes_{f} \phi} && K \otimes_{R^{'}} \ K. } \]

\noindent We then require that the following diagrams commute

\[  \xymatrix{ H \prescript{}{L}{\otimes} \ H \ar[rr]^-{\phi\prescript{}{f}{\otimes} \ \phi} \ar[dd]_-{\mu^{H}_{L}} && K \prescript{}{L^{'}}{\otimes} \ K \ar[dd]^-{\mu^{K}_{L}} \\
&& \\
H \ar[rr]_-{\phi} && K } \hspace{.25in}
\xymatrix{ H \otimes_{R} \ H \ar[rr]^-{\phi\ \otimes_{f} \phi} \ar[dd]_-{\mu^{H}_{R}} && K \otimes_{R^{'}} \ K \ar[dd]^-{\mu^{K}_{R}} \\
&& \\
H \ar[rr]_-{\phi} && K} \]

\item[(c)] Also by of condition $(a)$, the $k$-bimodule maps $\phi\prescript{}{f}{\otimes} \ \phi$ and $\phi\ \otimes_{f} \phi$ of condition $(b)$ further define $k$-bimodule maps

\[  \xymatrix{ H \prescript{}{L}{\times} \ H \ar[rr]^-{\phi\prescript{}{f}{\times} \ \phi} && K \prescript{}{L^{'}}{\times} \ K, } \hspace{.25in} \xymatrix{ H \times_{R} \ H \ar[rr]^-{\phi\ \times_{f} \phi} && K \times_{R^{'}} \ K. } \]

\noindent We then require that the following diagrams commute.

\[  \xymatrix{ H \ar[rr]^-{\phi} \ar[dd]_-{\Delta^{H}_{L}} && K \ar[dd]^-{\Delta^{K}_{L}} \\
&& \\
H \prescript{}{L}{\times} \ H \ar[rr]_-{\phi\prescript{}{f}{\times} \ \phi} && K \prescript{}{L^{'}}{\times} \ K } \hspace{.25in}
\xymatrix{ H \ar[rr]^-{\phi} \ar[dd]_-{\Delta^{H}_{R}} && K \ar[dd]^-{\Delta^{K}_{R}} \\
&& \\
H \times_{R} \ H \ar[rr]_-{\phi\ \times_{f} \phi} && K \times_{R^{'}} \ K } \]

\item[(d)] $\phi$ intertwines the antipodes of $\mathcal{H}$ and $\mathcal{K}$, i.e. $\phi\circ S_{H}=S_{K}\circ\phi$.

\end{enumerate}
\end{defn}

\begin{rem}\label{R2.3}
\begin{enumerate}
\item[]

\item[(1)] For a $k$-algebra $R$, let us denote by $HALG^{alg}(R)$ the category whose objects are Hopf algebroids over $R$ and morphisms are algebraic morphisms. For a fixed $k$, let us denote by $HALG^{geom}(k)$ the category whose objects are Hopf algebroids over $k$-algebras and morphisms are geometric morphisms. The existence of these two naturally defined categories reflect the fact that Hopf algebroids are generalization of both Hopf algebras and groupoids.

\item[(2)] Equip $R^{e}$ with the Hopf algebroid structure defined in example 5 of the next section. Let $(H_{L},H_{R},S)$ be a Hopf algebroid over $R$. Then the unit maps $\eta_{L},\eta_{R}$ together with the identity map on $R$ define geometric morphisms $(id,\eta_{L}):R^{e}\longrightarrow \mathcal{H}$ and $(id,\eta_{R}):R^{e}\longrightarrow \mathcal{H}$.
\end{enumerate}
\end{rem}

\subsection{Examples}\label{S2.2}

\begin{exa}\label{E2.1} \textbf{Hopf algebras}. A Hopf algebra $H$ over the commutative unital ring $k$ gives an example of a Hopf algebroid. Here, we take $R=L=k$ as $k$-algebras, take $s_{L}=t_{L}=s_{R}=t_{R}=\eta$ to be the source and target maps, set $\epsilon_{L}=\epsilon_{R}=\epsilon$ to be the counits, and $\Delta_{L}=\Delta_{R}=\Delta$ to be the coproducts.
\end{exa}

\begin{exa}\label{E2.2} \textbf{Coupled Hopf algebras}. It might be tempting to think that Hopf algebroids for which $R=L=k$ must be Hopf algebras. This is not entirely the case. We will give a general set of examples for which this is not true. Two Hopf algebra structures $H_{1}=(H,m,\eta,\Delta_{1},\epsilon_{1},S_{1})$ and $H_{2}=(H,m,\eta,\Delta_{2},\epsilon_{2},S_{2})$ over the same $k$-algebra $H$ are said to be \textit{coupled} if

\begin{enumerate}
\item[(a)] there exists a $k$-module map $C:H\longrightarrow H$, called the \textit{coupling map} such that

\[ \xymatrix{
& & H\otimes H \ar[r(3)]^-{C\otimes id} & & & & H\otimes H \ar[rd(1.6)]^-{m} & & \\
& & & & & & & & \\
H \ar[ru(1.8)]^-{\Delta_{1}} \ar[rd(1.8)]_-{\Delta_{2}} \ar@<1ex>[rrrr]^-{\epsilon_{2}} \ar@<-1ex>[rrrr]_-{\epsilon_{1}} & & & & k \ar[rrrr]^-{\eta} & & & & H \\
& & & & & & & & \\
& & H\otimes H \ar[r(3)]_-{id\otimes C} & & & & H\otimes H \ar[ru(1.6)]_-{m} & & \\
} \]

\noindent commutes, and

\item[(b)] the coproducts $\Delta_{1}$ and $\Delta_{2}$ in $H$ commutes.
\end{enumerate}

\noindent Coupled Hopf algebras give rise to Hopf algebroids over $k$. The left $k$-bialgebroid is the underlying bialgerba of $H_{1}$ while the right $k$-bialgebroid is the underlying bialgebra of $H_{2}$. The coupling map plays the role of the antipode.

Let us give examples of coupled Hopf algebras. Connes and Moscovici constructed \textit{twisted} antipodes in \cite{cm}. Let us show that such a twisted antipode is a coupling map for some coupled Hopf algebras. Let $H=(H,m,1,\Delta,\epsilon,S)$ be a Hopf algebra. Take $H_{1}=H$ as Hopf algebras. Let $\sigma:H\longrightarrow k$ be a character. Define $\Delta_{2}:H\longrightarrow H\otimes H$ by $h\mapsto h_{(1)}\otimes\sigma(S(h_{(2)}))h_{(3)}$. Take $\epsilon_{2}=\sigma$. Define $S_{2}:H\longrightarrow H$ by $h\mapsto \sigma(h_{(1)})S(h_{(2)})\sigma(h_{(3)})$. Note the Sweedler-legs of $h$ appearing in the definition of $S_{2}$ is the one provided by $\Delta$ and not by $\Delta_{2}$. Then, $H_{2}=(H,m,1,\Delta_{2},\epsilon_{2},S_{2})$ is a Hopf algebra coupled with $H_{1}$ by the coupling map $S^{\sigma}:H\longrightarrow H$ defined by $h\mapsto \sigma(h_{(1)})S(h_{(2)})$.
\end{exa}

\begin{exa}\label{E2.3} \textbf{Groupoid algebras}. Given a small groupoid $\mathcal{G}$ with finitely many objects and a commutative unital ring $k$, we can construct what is called the groupoid algebra of $\mathcal{G}$ over $k$, denoted by $k\mathcal{G}$. For such a groupoid $\mathcal{G}$, let us denote by $\mathcal{G}^{(0)}$ its set of objects, $\mathcal{G}^{(1)}$ its set of morphisms, $s,t:\mathcal{G}^{(1)}\longrightarrow \mathcal{G}^{(0)}$ the source and target maps, $\iota:\mathcal{G}^{(0)}\longrightarrow \mathcal{G}^{(1)}$ the unit map, $\nu:\mathcal{G}^{(1)}\longrightarrow \mathcal{G}^{(1)}$ the inversion map, $\mathcal{G}^{(2)}=\mathcal{G}^{(1)}\prescript{}{t}\times_{s} \ \mathcal{G}^{(1)}$ the set of composable pairs of morphisms, and $m:\mathcal{G}^{(2)}\longrightarrow \mathcal{G}^{(1)}$ the partial composition. The groupoid algebra $k\mathcal{G}$ is the $k$-algebra generated by $\mathcal{G}^{(1)}$ subject to the relation

\[ ff^{'}= 
\begin{cases}
    f\circ f^{'},& \text{if } f,f^{'} \ \text{are composable}\\
		& \\
    0,              & \text{otherwise}
\end{cases}
\]

\noindent for $f,f^{'}\in \mathcal{G}^{(1)}$. The groupoid algebra $k\mathcal{G}$ is a Hopf algebroid as folows. The base algebras $R$ and $L$ are both equal to $k\mathcal{G}^{(0)}$ and the two bialgebroids $H_{R}$ and $H_{L}$ are isomorphic as bialgebroids with underlying $k$-module $k\mathcal{G}^{(1)}$. The partial groupoid composition $m$ dualizes and extends to a multiplication $m:k\mathcal{G}^{(1)}\otimes k\mathcal{G}^{(1)}\longrightarrow k\mathcal{G}^{(1)}$ which then factors through the canonical surjection $k\mathcal{G}^{(1)}\otimes k\mathcal{G}^{(1)}\longrightarrow k\mathcal{G}^{(1)}\otimes_{k\mathcal{G}^{(0)}} k\mathcal{G}^{(1)}$ to give the product $k\mathcal{G}^{(1)}\otimes_{k\mathcal{G}^{(0)}} k\mathcal{G}^{(1)}\longrightarrow k\mathcal{G}^{(1)}$. The source and target maps $s,t$ of the groupoid give the source and target maps $s,t:k\mathcal{G}^{(0)}\longrightarrow k\mathcal{G}^{(1)}$, respectively. The unit map gives the counit map $\epsilon:k\mathcal{G}^{(1)}\longrightarrow k\mathcal{G}^{(0)}$. Finally, the inversion map gives the antipode map $S:k\mathcal{G}^{(1)}\longrightarrow k\mathcal{G}^{(1)}$. Note that the underlying bimodule structures of the right and the left bialgerboid is related by the antipode map.
\end{exa}

\begin{exa}\label{E2.4} \textbf{Weak Hopf algebras}. Another structure that generalize Hopf algebras, called weak Hopf algebras, also are Hopf algebroids. Explicitly, a weak Hopf algebra $H$ over a commutative unital ring $k$ is a unitary associative algebra together with $k$-linear maps $\Delta:H\longrightarrow H\otimes H$ (weak coproduct), $\epsilon:H\longrightarrow k$ (weak counit) and $S:H\longrightarrow H$ (weak antipode) satisfying the following axioms:

\begin{enumerate}
\item[(i)] $\Delta$ is multiplicative, coassociative, and weak-unital, i.e.
\[(\Delta(1)\otimes 1)(1\otimes \Delta(1))=\Delta^{(2)}(1)=(1\otimes\Delta(1))(\Delta(1)\otimes 1),\]

\item[(iii)] $\epsilon$ is counital, and weak-multiplicative, i.e. for any $x,y,z\in H$
\[ \epsilon(xy_{(1)})\epsilon(y_{(2)}z)=\epsilon(xyz)=\epsilon(xy_{(2)})\epsilon(y_{(1)}z),\]

\item[(v)] for any $h\in H$, $S(h_{(1)})h_{(2)}S(h_{(3)})=S(h)$ and
\[ h_{(1)}S(h_{(2)})=\epsilon(1_{(1)}h)1_{(2)}, \hspace{.5in} S(h_{(1)})h_{(2)}=1_{(1)}\epsilon(h1_{(2)}) \]
\end{enumerate}

Let us sketch a proof why a weak Hopf algebra $H$ is a Hopf algebroid. Consider the maps $p_{R}:H\longrightarrow H$, $h\mapsto 1_{(1)}\epsilon(h1_{(2)})$ and $p_{L}:H\longrightarrow H$, $h\mapsto \epsilon(1_{(1)}h)1_{(2)}$. By $k$-linearity and weak-multiplicativity of $\epsilon$, $p_{R}$ and $p_{L}$ are idempotents.

Multiplicativity and coassiociativity of $\Delta$ and counitality of $\epsilon$ implies that for any $h\in H$,

\[ h_{(1)}\otimes p_{L}(h_{(2)})=1_{(1)}h\otimes 1_{(2)} \hspace{.5in} p_{R}(h_{(1)})\otimes h_{(2)}=1_{(1)}\otimes h1_{(2)}.\]

\noindent Now, using these relations and coassiociativity of $\Delta$ we get

\[ 1_{(1)}1_{(1')}\otimes 1_{(2)} \otimes 1_{(2')} = 1_{(1')(1)}\otimes p_{L}(1_{(1')(2)})\otimes 1_{(2')}=1_{(1)}\otimes p_{L}(1_{(2)})\otimes 1_{(3)}\]

\[ 1_{(1)}\otimes 1_{(1')}\otimes 1_{(2)}1_{(2')} = 1_{(1)}\otimes p_{L}(1_{(2)(1)})\otimes 1_{(2)(2)}=1_{(1)(1)}\otimes p_{L}(1_{(1)(2)})\otimes 1_{(2)}\]

\noindent Thus, the first tensor factor of the left-hand side of the first equation above is in the image of $p_{R}$. Similarly, the last tensor factor of the left-hand side of the second equation above is in the image of $p_{L}$. Clearly, $p_{R}(1)=p_{L}(1)=1$. Hence, the images of $p_{R}$ and $p_{L}$ are unitary subalgebras of $H$. Denote these subalgebras by $R$ and $L$, respectively. By the weak-unitality of $\Delta$ we see that these subalgebras are commuting subalgebras of $H$.

Taking the source map $s$ as the inclusion $R\longrightarrow H$ and the target map as $t:R^{op}\longrightarrow H$, $r\mapsto\epsilon(r1_{(1)})1_{(2)}$ equips $H$ with an $R\otimes_{k}R^{op}$-ring structure. Taking $\epsilon_{R}=p_{R}$ and $\Delta_{R}$ as the composition

\[ \xymatrix{H \ar[rr]^-{\Delta}  & & H\otimes_{k}H \ar@{->>}[rr] & & H\otimes_{R}H }\]

\noindent equips $H$ with an $R$-coring structure $(H,\Delta_{R},\epsilon_{R})$. The ring and coring structures just constructed gives $H$ a structure of right $R$-bialgebroid $H_{R}$.

Using $R^{op}$ in place of $R$ in the above construction, we get a left $R^{op}$-bialgebroid $H_{R^{op}}$. Together with the right $R$-bialgebroid constructed and the existing weak antipode $S$, we get a Hopf algebroid $(H_{R^{op}},H_{R},S)$.
\end{exa}

\subsection{Representation of Hopf algebroids}\label{S2.3}

In this section, we will look at representations of Hopf algebroids. Towards the end of the section, we will look at the descent theoretic aspect of a special class of modules over Hopf algebroids, the so called relative Hopf modules. Let $\mathcal{H}=(H_{L},H_{R},S)$ be a Hopf algebroid with underlying $k$-module $H$. $H$ carries both a left $L$-module sctructure and a left $R$-module structure via the maps $s_{L}$ and $t_{R}$, respectively. A \textit{right} $\mathcal{H}$-\textit{comodule} $M$ is a right $L$-module and a right $R$-module together with a right $H_{R}$-coaction $\rho_{R}:M\longrightarrow M\otimes _{R}H$ and a right $H_{L}$-coaction $\rho_{L}:M\longrightarrow M\otimes_{L}H$ such that $\rho_{R}$ is an $H_{L}$-comodule map and $\rho_{L}$ is an $H_{R}$-comodule map.

For the coaction $\rho_{R}$, let us use the following Sweedler notation:

\[ \rho_{R}(m) = m^{[0]}\tens{R} m^{[1]} \]

\noindent and for the coaction $\rho_{L}$, let us use the following Sweedler notation:

\[ \rho_{L}(m) = m_{[0]} \tens{L} m_{[1]}. \]

\noindent With these notations, the conditions above explicitly means that for all $m\in M$, $l\in L$ and $r\in R$ we have

\[ (m\cdot l)^{[0]}\tens{R}(m\cdot l)^{[1]}=\rho_{R}(m\cdot l)=m^{[0]}\tens{R} t_{L}(l)m^{[1]} \]

\[ (m\cdot r)_{[0]}\tens{L}(m\cdot r)_{[1]}=\rho_{L}(m\cdot r)=m_{[0]}\tens{L} m_{[1]}s_{R}(r). \]

\noindent We further require that the two coactions satify the following commutative diagrams

\begin{equation}\label{eq2.1}
\xymatrix{
M \ar[d(1.7)]_-{\rho_{R}} \ar[r(1.7)]^-{\rho_{L}}
& & M\tens{L} H \ar[d(1.6)]^-{\rho_{R}\tens{L} id}\\
& & \\
M\tens{R} H \ar[r(1.6)]_-{id \tens{R}\Delta_{L}}
& & M\tens{R} H\tens{L} H } \hspace{0.75in}
\xymatrix{
M \ar[d(1.7)]_-{\rho_{L}} \ar[r(1.7)]^-{\rho_{R}}
& & M\tens{R} H \ar[d(1.6)]^-{\rho_{L}\tens{R} id}\\
& & \\
M\tens{L} H \ar[r(1.6)]_-{id \tens{L}\Delta_{R}}
& & M\tens{L} H\tens{R} H }
\end{equation}

We will denote by $\mathcal{M}^{\mathcal{H}}$ the category of right $\mathcal{H}$-comodules. Symmetrically, we can define left $\mathcal{H}$-comodules and we denote the category of a such by $^{\mathcal{H}}\mathcal{M}$.

Comodules over Hopf algebroids are comodules over the constituent bialgebroids. Thus, one can speak of two different coinvariants, one for each bialgebroid. For a given right $\mathcal{H}$-comodule $M$, they are defined as follows:

\[ M^{co \ H_{R}} = \left\{m\in M\left| \ \rho_{R}(m)=m\tens{R} 1\right.\right\}, \]

\[ M^{co \ H_{L}} = \left\{m\in M\left| \ \rho_{L}(m)=m\tens{L} 1\right.\right\}. \]

\noindent In the general case, we have $M^{co \ H_{R}}\subseteq M^{co \ H_{L}}$. But in our case, where we assume $S$ is bijective these two spaces coincide. This will be important in the formulation of Galois theory for Hopf algebroids. To see that these coinvariants coincide, consider the following map

\[ \Phi_{M}:M\tens{R}H\longrightarrow M\tens{L}H \]
\[ m\tens{R} h \mapsto \rho_{L}(m)\cdot S(h) \]

\noindent Here, $H$ acts on the right of $M\otimes_{L}H$ through the second factor. If $m\in M^{co \ H_{R}}$, then we have

\begin{eqnarray}
\nonumber\label{} \rho_{L}(m)&=&\rho_{L}(m)\cdot S(h) = \Phi_{M}(m\tens{R}1) = \Phi_{M}(\rho_{R}(m))\\
\nonumber &=& \Phi_{M}(m^{[0]}\tens{R} m^{[1]}) = \rho_{L}(m^{[0]})\cdot S(m^{[1]})\\
\nonumber &=& (m^{[0]}_{[0]}\tens{L} m^{[0]}_{[1]})\cdot S(m^{[1]}) = m^{[0]}_{[0]}\tens{L} m^{[0]}_{[1]}S(m^{[1]})\\
\nonumber &=& m_{[0]}\tens{L} m^{[0]}_{[1]}S(m^{[1]}_{[1]})= m_{[0]}\tens{L} s_{L}(\epsilon_{L}(m_{[1]}))\\
\nonumber &=& m_{[0]}s_{L}(\epsilon_{L}(m_{[1]}))\tens{L} 1 = m\tens{L} 1\\
\nonumber
\end{eqnarray}

\noindent This shows the inclusion $M^{co \ H_{R}}\subseteq M^{co \ H_{L}}$. To show the other inclusion, one can run the same computation but using the inverse of $\Phi_{M}$ which is the following map

\[ \Phi_{M}^{-1}:M\tens{L}H\longrightarrow M\tens{R}H \]
\[ m\tens{L}h\mapsto S^{-1}(h)\cdot\rho_{R}(m). \]

\noindent In this case, we can simply write $M^{co \ \mathcal{H}}$ for $M^{co \ H_{R}}=M^{co \ H_{L}}$ and refer to it as the $\mathcal{H}$-coinvariants of $M$ instead of distinguishing the $H_{R}$- from the $H_{L}$-coinvariants, unless it is necessary to do so.

Let us now discuss monoid objects in $\mathcal{M}^{\mathcal{H}}$. They are called $\mathcal{H}$-comodule algebras. A right $\mathcal{H}$-\textit{comodule algebra} is an $R$-ring $(M,\mu,\eta)$ such that $M$ is a right $\mathcal{H}$-comodule and $\eta:R\longrightarrow M$ and $\mu:M\otimes_{R}M\longrightarrow M$ are $\mathcal{H}$-comodule maps. Using Sweedler notation for coactions, this explicitly means that for any $m,n\in M$ we have

\begin{equation}\label{eq2.2}
(mn)^{[0]}\tens{R}(mn)^{[1]}=\rho_{R}(mn)=m^{[0]}n^{[0]}\tens{R} m^{[1]}n^{[1]},
\end{equation}

\begin{equation}\label{eq2.3}
(mn)_{[0]}\tens{L}(mn)_{[1]}=\rho_{L}(mn)=m_{[0]}n_{[0]}\tens{L} m_{[1]}n_{[1]},
\end{equation}

\begin{equation}\label{eq2.4}
1_{M}^{[0]}\tens{R}1_{M}^{[1]}=\rho_{R}(1_{M})=1_{M}\tens{R} 1_{H},
\end{equation}

\begin{equation}\label{eq2.5}
(1_{M})_{[0]}\tens{L}(1_{M})_{[1]}=\rho_{L}(1_{M})=1_{M}\tens{L} 1_{H}.
\end{equation}

Let $\mathcal{H}=(H_{L},H_{R},S)$ be a Hopf algebroid with underlying $k$-module $H$. A $k$-algebra extension $A\subseteq B$ is said to be (\textit{right}) $H_{R}$-\textit{Galois} if $B$ is a right $H_{R}$-comodule algebra with $B^{co \ H_{R}}=A$ and the map

\[ \xymatrix{B\tens{A}B \ar[rr]^-{\mathfrak{gal}_{R}} & & B\tens{R}H } \]
\[ a\tens{A}b\longmapsto ab^{[0]}\tens{R}b^{[1]}\]

\noindent is a bijection. The map $\mathfrak{gal}_{R}$ is called the Galois map associated to the bialgebroid extension $A\subseteq B$. Symmetrically, the extension $A\subseteq B$ is (\textit{right}) $H_{L}$-\textit{Galois} if $B$ is a right $H_{L}$-comodule algebra with $B^{co \ H_{L}}=A$ and the map

\[ \xymatrix{B\tens{A}B \ar[rr]^-{\mathfrak{gal}_{L}} & & B\tens{L}H } \]
\[ a\tens{A}b\longmapsto a_{[0]}b\tens{L}a_{[1]}\]

\noindent is a bijection. We say that a $k$-algebra extension $A\subseteq B$ is $\mathcal{H}$-\textit{Galois} if it is both $H_{R}$-Galois and $H_{L}$-Galois. It is not known in general if the bijectivity of $\mathfrak{gal}_{R}$ and $\mathfrak{gal}_{L}$ are equivalent. However, if the antipode $S$ is bijective (which is part of our standing assumption) then $\mathfrak{gal}_{R}$ is bijective if and only if $\mathfrak{gal}_{L}$. To see this, note that $\mathfrak{gal}_{L}=\Phi_{B} \circ \mathfrak{gal}_{R}$ where $\Phi_{B}$ is the map defined in the previous section for $M=B$. Since $S$ is bijective, $\Phi_{B}$ is an isomorphism which gives the desired equivalence of bijectivity of $\mathfrak{gal}_{R}$ and $\mathfrak{gal}_{L}$. Thus, the extension $A\subseteq B$ is $\mathcal{H}$-Galois if it is a bialgebroid Galois extension for any of its constituent bialgebroids.


\section{Hopf categories} \label{S3.0}

\subsection{Definitions and properties}\label{S3.1}

Batista et al. \cite{bcv} introduced the notion of a Hopf category over an arbitrary strict braided monoidal $\mathcal{V}$. In this section, we will introduce its topological version. For this purpose, we specialize $\mathcal{V}$ as the category of complex vector spaces whose braiding is the usual flip of tensor factors. Also, we will assume that the underlying categories of such Hopf categories are small. We will be primarily interested with \textit{finite-type} $\mathcal{V}$-enriched categories, by which we mean the hom-sets are finite-dimensional vector spaces. Before giving the definition of a Hopf category, let is introduce some notation first. For two $\mathcal{V}$-enriched categories $\mathscr{A}$ and $\mathscr{B}$ with the same set of objects $X$, we define $\mathscr{A}\otimes_{X}\mathscr{B}$ to the the $\mathcal{V}$-enriched category with $X$ as the set of objects and for $x,y\in X$, the hom-set of arrows from $x$ to $y$ is the vector space

\begin{equation}\label{eq3.1}
\left(\mathscr{A}\otimes_{X}\mathscr{B}\right)_{x,y}:=\mathscr{A}_{x,y} \otimes \mathscr{B}_{x,y}. \end{equation}

\noindent We call $\mathscr{A}\otimes_{X}\mathscr{B}$ the \textit{tensor product} of $\mathscr{A}$ and $\mathscr{B}$. With this $\otimes_{X}$, the category of $\mathcal{V}$-enriched categories over $X$ becomes a strict monoidal category whose monoidal unit, denoted by $\mathbb{1}^{X}$, is the category over $X$ such that for any $x,y\in X$ we have $\mathbb{1}^{X}_{x,y}=\mathbb{C}$.

\begin{defn}\label{D3.1}
A \textit{Hopf category} $\mathscr{H}$ over $X$ is a $\mathcal{V}$-enriched category satisfying the following conditions.

\begin{enumerate}

\item[(a)] There are functors

\[ \xymatrix{\mathscr{H} \ar[rr]^-{\Delta} && \mathscr{H}\tens{X} \mathscr{H}}, \hspace{.5in} \xymatrix{\mathscr{H} \ar[rr]^-{\epsilon} && \mathbb{1}^{X}} \]

\noindent called the \textit{coproduct} and \textit{counit}, respectively, such that $\Delta$ is \textit{coassociative} and \textit{counital} with respect to $\epsilon$, i.e. the diagram of functors

\[ \xymatrix{ \mathscr{H} \ar[rr]^-{\Delta} \ar[dd]_-{\Delta} && \mathscr{H}\tens{X} \mathscr{H} \ar[dd]|-{id\otimes_{X}\Delta} \\
&& \\
\mathscr{H}\tens{X} \mathscr{H} \ar[rr]_-{\Delta\otimes_{X} id} && \mathscr{H}\tens{X} \mathscr{H}\tens{X} \mathscr{H} \\ }
\hspace{.25in}
\xymatrix{ \mathscr{H} \ar@{=}[rr] \ar@{=}[dd] \ar[rrdd]^-{\Delta} && \mathscr{H}\tens{X} \mathbb{1}^{X} \\
&& \\
\mathbb{1}^{X}\tens{X} \mathscr{H} && \mathscr{H}\tens{X} \mathscr{H} \ar[ll]^-{\epsilon\otimes_{X} id} \ar[uu]_-{is\otimes_{X}\epsilon} \\ }
 \]

\noindent commute.

\item[(b)] There is a functor $S:\mathscr{H}\longrightarrow \mathscr{H}^{op}$, called the \textit{antipode}, satisfying

\[ \xymatrix{
&& \mathscr{H}\tens{X} \mathscr{H} \ar[rr]^-{S\otimes_{X} id} & & \mathscr{H}^{op}\tens{X} \mathscr{H} \ar[rrd]^-{\circ} && \\
\mathscr{H} \ar[rru]^-{\Delta} \ar[rrd]_-{\Delta} \ar[rrr]^-{\epsilon} & & & \mathbb{1}_{X} \ar[rrr]^-{\eta} & & & \mathscr{H} \\
&& \mathscr{H}\tens{X} \mathscr{H} \ar[rr]_-{id\otimes_{X} S} & & \mathscr{H}\tens{X} \mathscr{H}^{op} \ar[rru]_-{\circ} && \\} \]

\noindent Here, $\circ$ denotes the bifunctor induced by the categorical composition in $\mathscr{H}$ and $\eta$ is the functor that send $1\in\mathbb{1}^{X}_{x,y}$ to the identity element of $\mathscr{H}_{x,y}$.
\end{enumerate}
\end{defn}

\begin{rem}\label{R3.1}
Functoriality of $\Delta$ and $\epsilon$ means that for any $x,y\in X$, we have linear maps

\[ \xymatrix{\mathscr{H}_{x,y} \ar[rr]^-{\Delta_{x,y}} && \mathscr{H}_{x,y}\otimes \mathscr{H}_{x,y}} \hspace{.5in} \xymatrix{\mathscr{H}_{x,y} \ar[rr]^-{\epsilon_{x,y}} && \mathbb{C}} \]

\noindent where $\Delta_{x,y}$ is coassociative and counital with respect to $\epsilon_{x,y}$ in the usual sense. This implies that $\mathscr{H}_{x,y}$ is a coalgebra. If we denote by $C{V}$ the category of coalgebras on $\mathcal{V}$, another way to package part $(a)$ of definition (\ref{D3.1}) is to say that $\mathscr{H}$ is enriched over $C(\mathcal{V})$.
\end{rem}

For the main results of this paper, we will be mostly interested with the case $X$ is a topological space. In such a case, it makes sense to reflect \textit{continuity} on the functors $\Delta$, $\epsilon$ and $S$ along with the categorical structure maps. This calls for the following definition.

\begin{defn}\label{D3.2}
Let $X$ be a topological space and let $\mathcal{O}_{X}$ be the sheaf of continuous complex-valued functions on $X$. A \textit{topological Hopf category} $\mathscr{H}$ over $X$ is a Hopf category together with a sheaf $H^{sh}$ over $X\times X$ (with the product topology) of $\mathcal{O}_{X}$-bimodules satisfying the following conditions.
\begin{enumerate}
\item[(a)] Denote by $\pi_{1},\pi_{2}:X\times X\longrightarrow X$ the projection onto the first and second factor, respectively. Over an open set $U\subseteq X\times X$, for any $\sigma\in H^{sh}(U)$, $f\in \mathcal{O}_{X}(\pi_{1}U)$ and $g\in \mathcal{O}_{X}(\pi_{2}U)$ we have

\[ \left(f\cdot\sigma\cdot g\right)(x,y)=f(x)\sigma(x,y)g(y) \]

\noindent for any $(x,y)\in U$.

\item[(b)] $\mathscr{H}_{x,y}$ is the fiber of $H^{sh}$ at $(x,y)\in X\times X$.

\item[(c)] $\circ$, $\eta$, $\Delta$, $\epsilon$ and $S$ are the induced maps on global sections of the following map of sheaves

\[ \xymatrix{H^{sh}\tens{\mathcal{O}_{X}} H^{sh} \ar[rr]^-{\circ^{sh}} && H^{sh} }, \hspace{.5in} \xymatrix{\mathcal{O}_{X} \ar[rr]^-{\eta^{sh}} && H^{sh} },  \]

\[ \xymatrix{H^{sh} \ar[rr]^-{\Delta^{sh}} && H^{sh}\tens{\mathcal{O}_{X}} H^{sh}}, \hspace{.5in} \xymatrix{H^{sh} \ar[rr]^-{\epsilon^{sh}} && \mathcal{O}_{X}},  \]

\[ \xymatrix{ H^{sh} \ar[rr]^-{S^{sh}} && \left(H^{sh}\right)^{op} } \]

\noindent respectively. Here, $\left(H^{sh}\right)^{op}$ is the pullback of the sheaf $H^{sh}$ along the map $X\times X\longrightarrow X\times X$ flipping the factors.
\end{enumerate}
\end{defn}

\begin{rem}\label{R3.2}
The bimodule tensor product $\otimes_{\mathcal{O}_{X}}$ used in part $(c)$ for $\Delta^{sh}$ of definition (\ref{D3.2}) is the tensor product of the appropriately modified $\mathcal{O}_{X}$-bimodule $H^{sh}$, one in which we have

\[ f\cdot \left(\sigma \tens{\mathcal{O}_{X}} \tau\right)\cdot g = \left(\sigma\cdot g\right) \tens{\mathcal{O}_{X}} \left(f\cdot\tau\right) \]

\noindent for any $f,g\in \mathcal{O}_{X}$. For the bimodule tensor product $\otimes_{\mathcal{O}_{X}}$ used for $\circ^{sh}$ is the one with

\[ f\cdot \left(\sigma \tens{\mathcal{O}_{X}} \tau\right)\cdot g = \left(f\cdot\sigma\right) \tens{\mathcal{O}_{X}} \left(\tau\cdot g\right) \]

\noindent for any $f,g\in\mathcal{O}_{X}$.
\end{rem}

The following, which will play an important role in our formulation of the main result, is the categorification of a coupled Hopf algebra.

\begin{defn}\label{D3.3}
A \textit{coupled Hopf category} $\mathscr{H}$ is a $\mathcal{V}$-enriched category with two $C(\mathcal{V})$-enrichments, denoted by $\mathscr{H}_{L}$ and $\mathscr{H}_{R}$, with coproducts $\Delta^{L},\Delta^{R}$ and counits $\epsilon^{L},\epsilon^{R}$, respectively; and a functor $S:\mathscr{H}\longrightarrow \mathscr{H}^{op}$, called the \textit{coupling} functor, such that the following conditions are satisfied:

\begin{enumerate}
\item[(a)] The following diagrams, indicating the \textit{coupling condition}, commute.

\[ \xymatrix{
& & \mathscr{H}\tens{X} \mathscr{H} \ar[rr]^-{S\otimes_{X} id} & & \mathscr{H}^{op}\tens{X} \mathscr{H} \ar[rrd]^-{\circ} & & \\
\mathscr{H} \ar[rru]^-{\Delta^{L}}  \ar[rrr]^-{\epsilon^{R}} & & & \mathbb{1}^{X} \ar[rrr]^-{\eta} & & & \mathscr{H} \\
\mathscr{H} \ar[rrd]_-{\Delta^{R}}  \ar[rrr]_-{\epsilon^{L}} & & & \mathbb{1}^{X} \ar[rrr]_-{\eta} & & & \mathscr{H} \\
& & \mathscr{H}\tens{X} \mathscr{H} \ar[rr]_-{id\otimes_{X} S} & & \mathscr{H}\tens{X} \mathscr{H}^{op} \ar[rru]_-{\circ} & & \\
}\]

\item[(b)] The coproducts $\Delta^{L}$ and $\Delta^{R}$ commute, i.e.

\[ \xymatrix{
\mathscr{H} \ar[dd]_-{\Delta^{L}} \ar[rr]^-{\Delta^{R}}
& & \mathscr{H}\tens{X} \mathscr{H} \ar[dd]^-{\Delta^{L}\tens{X} id}\\
& & \\
\mathscr{H}\tens{X} \mathscr{H} \ar[rr]_-{id \tens{X}\Delta^{R}}
& & \mathscr{H}\tens{X} \mathscr{H}\tens{X} \mathscr{H} }
\hspace{.25in}
\xymatrix{
\mathscr{H} \ar[dd]_-{\Delta^{R}} \ar[rr]^-{\Delta^{L}}
& & \mathscr{H}\tens{X} \mathscr{H} \ar[dd]^-{\Delta^{R}\tens{X} id}\\
& & \\
\mathscr{H}\tens{X} \mathscr{H} \ar[rr]_-{id \tens{X}\Delta^{L}}
& & \mathscr{H}\tens{X} \mathscr{H}\tens{X} \mathscr{H} } \]

\end{enumerate}
\end{defn}

\begin{rem}\label{R3.3}
\begin{enumerate}
\item[]
\item[(1)] Coupled Hopf categories are almost the categorification of coupled Hopf algebras. While the constituent bialgebras of a coupled Hopf algebra is a Hopf algebras in itself, the constituent categories $\mathscr{H}_{L}$ and $\mathscr{H}_{R}$ of a coupled Hopf category $\mathscr{H}$ need not be Hopf categories.

\item[(2)] Just like Hopf categories, we can also \textit{topologize} coupled Hopf categories. We can take definition (\ref{D3.2}): assert the existence of a sheaf $H^{sh}$ over $X\times X$ of $\mathcal{O}_{X}$-bimodules, take conditions $(a)$ and $(b)$ as they are, and replace condition $(c)$ by
\begin{enumerate}
\item[(c')] $\Delta^{L}$, $\Delta^{R}$, $\epsilon^{L}$, $\epsilon^{R}$ and $S$ are the induced maps on global sections of the following map of sheaves

\[ \xymatrix{H^{sh} \ar[rr]^-{(\Delta^{L})^{sh}} && H^{sh}\tens{\mathcal{O}_{X}} H^{sh}}, \hspace{.5in} \xymatrix{H^{sh} \ar[rr]^-{(\epsilon^{L})^{sh}} && \mathcal{O}_{X}},  \]

\[ \xymatrix{H^{sh} \ar[rr]^-{(\Delta^{R})^{sh}} && H^{sh}\tens{\mathcal{O}_{X}} H^{sh}}, \hspace{.5in} \xymatrix{H^{sh} \ar[rr]^-{(\epsilon^{R})^{sh}} && \mathcal{O}_{X}},  \]

\[ \xymatrix{ H^{sh} \ar[rr]^-{S^{sh}} && \left(H^{sh}\right)^{op} } \]

\noindent respectively, making the following diagram

\begin{flushleft}
\xymatrix{
& H^{sh}(U)\tens{\mathcal{O}_{X}(U)} H^{sh}(U) \ar[rr]^-{S_{U}\tens{\mathcal{O}_{X}(U)} id} & & H^{sh}(U)^{op}\tens{\mathcal{O}_{X}(U)} H^{sh}(U) \ar[rd]^-{\mu_{U}} & \\
H^{sh}(U) \ar[ru]^-{\left(\Delta^{L}\right)^{sh}(U)}  \ar[rr]^-{(\epsilon^{R})^{sh}(U)} & & \mathcal{O}_{X}(U) \ar[rr]^-{\eta_{U}} & & H^{sh}(\pi^{diag}_{2}U) \\
H^{sh}(U) \ar[rd]_-{(\Delta^{R})^{sh}(U)}  \ar[rr]_-{(\epsilon^{L})^{sh}(U)} & & \mathcal{O}_{X}(U) \ar[rr]_-{\eta_{U}} & & H^{sh}(\pi^{diag}_{1}U) \\
& H^{sh}(U)\tens{\mathcal{O}_{X}(U)} H^{sh}(U) \ar[rr]_-{id\tens{\mathcal{O}_{X}(U)} S_{U}} & & H^{sh}(U)\tens{\mathcal{O}_{X}(U)} (H^{sh})(U)^{op} \ar[ru]_-{\mu_{U}} & \\}
\end{flushleft}

\noindent commute for any $U\subseteq X\times X$. Here, $\mu_{U}$ and $\eta_{U}$ denote the maps induced by the composition and unit maps of $\mathscr{C}$. The maps $\pi_{1}^{diag}$ and $\pi_{2}^{diag}$ denote $X\times X\longrightarrow X\times X$, $(x,y)\mapsto(x,x)$ and $X\times X\longrightarrow X\times X$, $(x,y)\mapsto(y,y)$, respectively.
\end{enumerate}

\end{enumerate}
\end{rem}

\subsection{A good example of a Hopf category}\label{S3.2}

In this section, we will look at a very important example of a Hopf category. This example will also be an example of our main result. This is a special case of proposition 7.1 of \cite{bcv}. Consider a finite set $X$ whose elements are conveniently labelled as $1,2,...,n$. Equipped $X$ with the discrete topology. Consider the category $C$ whose set of objects is $X$ and define $C_{x,y}=\mathbb{C}$. The category $C$ is obviously a Hopf category. By proposition 7.1 of \cite{bcv}, $\mathcal{H}=\bigoplus_{x,y\in X}C_{x,y}$ is a weak Hopf algebra. Using the arguments in example 4 of section (\ref{S2.2}), $\mathcal{H}$ is a Hopf algebroid over $A=\mathbb{C}^{n}=\mathcal{O}_{X}(X)$.

The Hopf algebroid $\mathcal{H}$ has a more familiar form. It is isomorphic, as a Hopf algebroid, the algebra $M_{n}(\mathbb{C})$ over its diagonal $D_{n}=Diag_{n}(\mathbb{C})$. With the $D_{n}$-bimodule structure on $M_{n}(\mathbb{C})$ defined as

\[ P\cdot M\cdot Q := MPQ, \hspace{.5in} P,Q\in D_{n}, M\in M_{n}(\mathbb{C}), \]

\noindent the coproduct $\Delta_{R}$ and the counit $\epsilon_{R}$ are given as

\[ \Delta_{R}(E_{ij})=E_{ij}\otimes_{D_{n}}E_{ij}, \hspace{.35in} \epsilon_{R}(M)=\sum\limits_{i=1}^{n}E_{ii}\phi(ME_{ii}) \]

\noindent where $\phi$ is the linear functional defined by $\phi(E_{ij})=1$ for all $i,j\in X$. With the usual matrix multiplication and unit, $\Delta_{R}$ and $\epsilon_{R}$ constitutes a right $D_{n}$-bialgebroid structure on $M_{n}(\mathbb{C})$. For completeness, let us define the structure maps of the left $D_{n}$-bialgebroid structure of $M_{n}(\mathbb{C})$. Consider the $D_{n}$-bimodule structure on $M_{n}(\mathbb{C})$ defined as

\[ P\cdot M\cdot Q := PQM, \hspace{.5in} P,Q\in D_{n}, M\in M_{n}(\mathbb{C}). \]

\noindent The coproduct $\Delta_{L}$ and the counit $\epsilon_{L}$ are defined as

\[ \Delta_{L}(E_{ij})=E_{ij}\otimes_{D_{n}}E_{ij}, \hspace{.35in} \epsilon_{L}(M)=\sum\limits_{i=1}^{n}\phi(E_{ii}M)E_{ii} \]

\noindent where $\phi$ is the same linear functional used to defined $\epsilon_{R}$. The antipode $S$ of this Hopf algebroid is defined as $S(E_{ij})=E_{ji}$.

As a weak Hopf algebra, $\phi$ is the counit of $\mathcal{H}$. The coproducts $\Delta_{L}$ and $\Delta_{R}$ are the extension of the weak coproduct $\Delta$ to $M_{n}(\mathbb{C})\otimes_{D_{n}}M_{n}(\mathbb{C})$ relative to the $D_{n}$-bimodule structure used. As we will see in section (\ref{S4.0}), this is not a coincidence. This is in fact a special case of a more general result which we shall prove at the end of that section.

\subsection{Galois extensions of Hopf categories}\label{S3.3}

Formulation of Galois theory for Hopf category is straightforward. Recall that in the case of Hopf algebras, only the underlying bialgebra structure is relevant. In the coaction picture, the coalgebra is used to make sense of a coaction while the algebra structure is used to make sense of the Galois map. All these ingredients are already present in the case of a Hopf category. We will discuss the situation for topological Hopf categories. The case for Hopf categories follow almost immediately by dropping any manifestation of topology.

Before giving the definition of the categorical analogue of a comodule algebra, let us first discuss what a topological category is, at least for our purpose. A $\mathcal{V}$-enriched category $\mathscr{M}$ over a space $X$ is a \textit{topological category} if there is a sheaf $M^{sh}$ of $\mathcal{O}_{X}$-bimodules such that conditions $(a)$, $(b)$ and the relevant part of condition $(c)$ of definition (\ref{D3.2}) hold.

\begin{defn}\label{D3.4} Let $\mathscr{H}$ be a topological Hopf category with space of objects $X$, coproduct $\Delta$, counit $\epsilon$ and antipode $S$ with associated sheaf $H^{sh}$.

\begin{enumerate}
\item[(1)] A topological category $\mathscr{M}$ over $X$ enriched over $\mathcal{V}$, with associated sheaf $M^{sh}$, is a \textit{right} $\mathscr{H}$-\textit{comodule} if there is a functor $\rho:\mathscr{M}\longrightarrow \mathscr{M}\otimes_{X}\mathscr{H}$ such that the following conditions hold.

\begin{enumerate}
\item[(a)] $\rho$ is coassociative with respect to $\Delta$ and counital with respect to $\epsilon$, i.e. the diagrams of functors

\[ \xymatrix{ \mathscr{M} \ar[rr]^-{\rho} \ar[dd]_-{\rho} && \mathscr{M}\otimes_{X}\mathscr{H} \ar[dd]|-{id\otimes_{X}\Delta} \\
&&& \\
\mathscr{M}\otimes_{X}\mathscr{H} \ar[rr]_-{\rho\otimes_{X} id} && \mathscr{M}\otimes_{X}\mathscr{H}\otimes_{X}\mathscr{H} \\ }
\hspace{-.5in}
\xymatrix{ \mathscr{M} \ar@{=}[rr] \ar[rdd]_-{\rho} && \mathscr{M}\otimes_{X}\mathbb{1}^{X} \\
&& \\
& \mathscr{M} \otimes_{X} \mathscr{H} \ar[ruu]|-{id\otimes_{X}\epsilon} & \\} \]

\noindent commute, and

\item[(b)] the functor $\rho$ is the map induced by the map of sheaves $M^{sh}\longrightarrow M^{sh}\otimes_{\mathcal{O}_{X}}H^{sh}$ where the tensor product is the same as the first one we described in remark (\ref{R3.2}). A \textit{left} $\mathscr{H}$-\textit{comodule} can be symmetrically defined.
\end{enumerate}

\item[(2)] A \textit{morphism} $\mathscr{M}\stackrel{\phi}{\longrightarrow}\mathscr{N}$ of right $\mathscr{H}$-comodules is a functor that commutes with the right coactions, i.e. one which makes the following diagram commute

\[ \xymatrix{ \mathscr{M} \ar[rrr]^-{\rho^{M}} \ar[dd]_-{\phi} &&& \mathscr{M}\otimes_{X}\mathscr{H} \ar[dd]|-{\phi\otimes_{X}id} \\
&&& \\
\mathscr{N} \ar[rrr]_-{\rho^{N}} &&& \mathscr{N}\otimes_{X}\mathscr{H}. \\ } \]

\noindent Here, $\rho^{M}$ and $\rho^{N}$ are the coactions of $\mathscr{H}$ on $\mathscr{M}$ and $\mathscr{N}$, respectively.

\item[(3)] A right $\mathscr{H}$-comodule $\mathscr{M}$ is a \textit{right} $\mathscr{H}$-\textit{comodule-category} if in addition, the composition map $\mathscr{M}\otimes_{X}\mathscr{M}\stackrel{\circ}{\longrightarrow} \mathscr{M}$ is a map of right $\mathscr{H}$-comodules, where $\mathscr{M}\otimes_{X}\mathscr{M}$ is equipped with the diagonal coaction.

\item[(4)] The \textit{coinvariants} of a right $\mathscr{H}$-comodule-category $\mathscr{M}$ is the subcategory $\mathscr{M}^{co\ \mathscr{H}}$ whose space of objects is $X$ and whose hom-sets are defined as

\[ \left(\mathscr{M}^{co\ \mathscr{C}}\right)_{x,y}:=\left\{ \alpha\in \mathscr{M}_{x,y}|\rho(\alpha)=\alpha\otimes id_{y} \right\} \]

\noindent for any $x,y\in X$.
\end{enumerate}
\end{defn}

\begin{rem}\label{R3.4}
A Hopf category is the categorification of a Hopf algebra with categorical composition corresponding to the algebra product. A right $\mathscr{H}$-comodule $\mathscr{M}$ is in particular a category, it already has a composition. This means that we only need to impose requirement $(3)$ in definition (\ref{D3.4}) to get a categorification of the notion of a comodule-algebra. In the classical set-up, one has to require the existence of a product and assert its compatibility with the comodule structures.
\end{rem}

In the set-up of Hopf-Galois theory with respect to Hopf algebras, there is a well-understood notion for extensions of $k$-algebras $A\subseteq B$ to be $H$-Galois for a Hopf algebra $H$ even if $A\neq k$. This is because $B\otimes_{A}B$ makes sense as a $k$-module. All that is left to do is require $A=B^{co\ H}$ and that the map $B\otimes_{A}B\longrightarrow B\otimes H$, $a\otimes b\mapsto (a\otimes 1)\rho(b)$ is bijective. On the other hand, in the situation of a Hopf category $\mathscr{H}$ and extensions of comodule-categories $\mathscr{A}\subseteq\mathscr{M}$ with $\mathscr{A}=(\mathscr{M})^{co\ \mathscr{H}}$, we can only make sense of the product $\mathscr{M}\otimes_{\mathscr{A}}\mathscr{M}$ in the case $\mathscr{A}$ is the subcategory of $\mathscr{M}$ whose hom-sets $\mathscr{A}_{x,y}$ are all zero except when $x=y$, in which case $\mathscr{A}_{x,x}=\mathbb{C}$. In this case, we identify $\mathscr{M}\otimes_{\mathscr{A}}\mathscr{M}$ with $\mathscr{M}\otimes_{X}\mathscr{M}$. Let us call such a category the \textit{trivial linear category} over $X$, and denote by $I_{X}$. There might be a way to consider Galois extensions by Hopf categories in which the subcategory of coinvariants is strictly larger than $I_{X}$, but at present it is not clear to the author how to make sense of it. Fortunately, for our purpose of proving theorem (\ref{T5.1}) it is enough to have $I_{X}$ as the subcategory of coinvariants.

\begin{defn}\label{D3.5}
A right $\mathscr{H}$-comodule-category $\mathscr{M}$ is a $\mathscr{H}$-\textit{Galois extension} of $I_{X}$ provided
\begin{enumerate}
\item[(a)] $\mathscr{M}^{co\ \mathscr{H}}=I_{X}$, and

\item[(b)] the functor

\[ \mathscr{M}\otimes_{X}\mathscr{M} \stackrel{\mathfrak{gal}}{\longrightarrow} \mathscr{M}\otimes_{X}\mathscr{H}, \]
\[ \alpha \otimes \beta \mapsto \left(\alpha\circ\beta_{[0]}\right)\otimes \beta_{[1]}\]

\noindent called the \textit{Galois morphism}, is fully faithful.
\end{enumerate}
\end{defn}

\begin{rem}\label{R3.5}
\begin{enumerate}
\item[]
\item[(1)] We are using Sweedler notation for the legs of the coaction

\[ \rho:\mathscr{M}\longrightarrow \mathscr{M}\otimes_{X}\mathscr{H}. \]

\noindent In other words, for any $x,y\in X$ and $\alpha\in\mathscr{M}_{x,y}$, we have $\rho(\alpha)=\alpha_{[0]}\otimes\alpha_{[1]}$, where $\alpha_{[0]}\in\mathscr{M}_{x,z}$ and $\alpha_{[1]}\in\mathscr{H}_{z,y}$ for some $z\in X$. This, in particular, tells us that the map $\mathfrak{gal}$ above make sense.
\item[(2)] Galois extension by a coupled Hopf category $\mathscr{H}=(\mathscr{H}_{L},\mathscr{H}_{R},S)$ means simultaneous Galois extensions of the constituent $C(\mathcal{V})$-enriched categories $\mathscr{H}_{L}$ and $\mathscr{H}_{R}$.
\end{enumerate}
\end{rem}


\section{The category associated to a Hopf algebroid} \label{S4.0}

In this section, we will consider Hopf algebroids $\mathcal{H}$ over a commutative unital $C^{*}$-algebra $A$. We will restrict to the case where $\mathcal{H}$ is finitely-generated and projective as a left and a right $A$-module. With the underlying assumption that the antipode is bijective, by \cite{bohm-szlachanyi}, finitely-generated projectivity of any of the $A$-module structures of $\mathcal{H}$ coming from the source and target maps are all equivalent. Note that even though $A$ is commutative, its image under the source or the target map need not be central in $\mathcal{H}$. We will deal with this general situation and specialize in the case when we have centrality.

\subsection{Local eigenspace decomposition} \label{S4.1}

Let $\mathcal{H}=(H_{L},H_{R},S)$ be a Hopf algebroid over $A$, a commutative unital $C^{*}$-algebra. Assume that $H_{L}$ is finitely-generated and projective as a left- and a right-$A$-module via the source and target maps. With our standing assumption, $H_{R}$ has the same properties.

Let us first consider the left bialgebroid $H_{L}$. The Gelfand duality implies that $A\cong C(X)$ for some compact Hausdorff space $X$. The Serre-Swan theorem applied to the left $A$-module $H_{L}$ gives us a finite-rank vector bundle $E\stackrel{p}{\longrightarrow}X$ such that $H_{L}\cong\Gamma(X,E)$ as left modules, where the left $C(X)$-module structure on $\Gamma(X,E)$ is by pointwise multiplication, i.e. $(f\cdot\sigma)(x)=f(x)\sigma(x)$ for all $x\in X, f\in C(X)$ and $\sigma\in \Gamma(X,E)$. By the bimodule nature of $H_{L}$, the right $A$-module structure of $\Gamma(X,E)$ commutes with the left $A$-module which implies that we have a representation $C(X)\stackrel{\rho}{\longrightarrow} End(E)$ of $C(X)$ into the endomorphism bundle of $E\stackrel{p}{\longrightarrow}X$. Since $C(X)$ is abelian and $\rho$ is a $*$-morphism, $\rho(C(X))$ lands in a maximal abelian subalgebra $D(n)$ of $End(E)$.

Choose a finite collection of open sets $\left\{U_{i}|i=1,2,...,m\right\}$ that cover $X$ over which $E$ is trivializable. Choose a system of coordinates such that $E$ trivial over each $U_{i}$, i.e. $\left.E\right|_{U_{i}}\cong U_{i}\times V$, where $V$ is a finite-dimensional vector space.Choosing a basis $v_{1},v_{2},...,v_{n}\in V$ one has $End(\left.E\right|_{U_{i}})=C(U_{i},M_{n}(\mathbb{C}))$ where $n$ is the rank of $E$. Commutativity of $C(X)$ implies that up to unitaries $V_{i}\in U(n)$, we have

\[ \xymatrix{ C(X) \ar[rr]^-{\rho} && C(U_{i}, Diag(n))} \]

\noindent where $Diag(n)$ denotes the subalgebra of diagonal matrices on $M_{n}(\mathbb{C})$ and

\[ V_{i}\cdot C(U_{i}, Diag(n)) \cdot V_{i}^{*}=\left.D(n)\right|_{U_{i}}. \]

\noindent For each $i=1,2,...,m$, choosing a set of central orthogonal idempotents $\left\{e_{j}|j=1,...,n\right\}$ gives $n$ projections $p^{i}_{j}$ given by the following composition

\[ \xymatrix{ C(X) \ar[rr]^-{\rho_{i}} && C(U_{i}, Diag(n))\cong \bigoplus\limits_{k=1}^{n}C(U_{i}) \ar[rr]^-{proj_{j}} && C(U_{i}) }  \]

\begin{figure}[ht!]
\centering
\includegraphics[width=135mm]{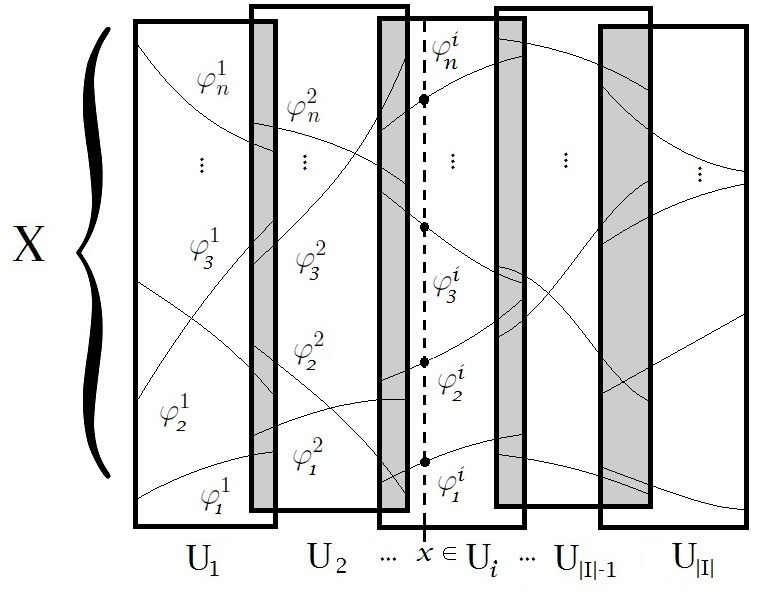}
\caption{Local eigenspace decomposition of $E$.}
\label{fg4.1}
\end{figure}

\noindent These projections are in particular continuous $C^{*}$-morphisms. Hence, they give, for each $i=1,2,...,m$, (possibly non-distinct) $n$ continuous injective maps $U_{i}\stackrel{\varphi^{i}_{j}}{\longrightarrow}X$, $j=1,...,n$. Geometrically, the situation is depicted figure (\ref{fg4.1}).

Let us describe the nature of the set $Z=\bigcup_{i,j}\varphi_{j}^{i}(U_{i})$ over the intersections $U_{\alpha}\cap U_{\beta}$. Over $U_{\alpha}\cap U_{\beta}\subseteq U_{\alpha}$ we get a unitary $V_{\alpha}$ which gives $n$ central orthogonal idempotents and up to ordering of such idempotents, one gets the sets $\varphi^{i}_{j}(U_{i})$. The union $\bigcup_{j}\varphi^{i}_{j}(U_{i})$ does not depend on the ordering of these idempotents. Thus, over $U_{\alpha}\cap U_{\beta}$ one gets unitaries $V_{\alpha}$ and $V_{\beta}$ which simultaneously diagonalize $\rho(C(X))$. Thus, we have

\[ \bigcup_{j}\varphi^{\alpha}_{j}(U_{\alpha}\cap U_{\beta}) = \bigcup_{j}\varphi^{\beta}_{j}(U_{\alpha}\cap U_{\beta}) \]

\noindent from which we get that

\[ \left(\bigcup_{j}\varphi^{\alpha}_{j}(U_{\alpha})\right)\cap\left(\bigcup_{j}\varphi^{\beta}_{j}(U_{\beta})\right)= \bigcup_{j}\varphi^{\alpha}_{j}(U_{\alpha}\cap U_{\beta}) \]

\noindent that is, the sets $\bigcup_{j}\varphi^{i}_{j}(U_{i})$ agree on the intersections.

A subset $T\subseteq X\times X$ is called \textit{transverse} if

\[ \left.proj_{1}\right|_{T}: X\times X\longrightarrow X, \hspace{.5in} \left.proj_{2}\right|_{T}: X\times X\longrightarrow X \]

\noindent are homeomorphisms, where $proj_{1}$ and $proj_{2}$ denotes the projection onto the first and second factor, respectively. In particular, $T$ is homeomorphic to $X$. Using the above argument, we have the following proposition.

\begin{prop}\label{P4.1}
For every $i=1,2,...,m$, $j=1,2,...,n$ the set $\varphi_{j}^{i}(U_{i})$ extends to a transverse subset of $X\times X$ completely contained in $Z$. In particular, $Z$ is the union of $n$ (possibly overlapping) transverse subsets of $X\times X$.
\end{prop}

\noindent This means that the curves in figure (\ref{fg4.1}) overlap.

\begin{rem}\label{R4.1}
Another way to see why the closed subset $Z\subset X\times X$ is the union of transverse subsets of $X\times X$ is by the fact the we can run the construction of the sets $\varphi^{i}_{j}(U_{i})$ described in the beginning of this section in a symmetric fashion, one for each factor of $X\times X$.
\end{rem}

The whole picture (\ref{fg4.1}) is a decomposition of $X\times X$ into $X\times U_{i}$, $i\in I$. The graphs of $\varphi^{i}_{j}$ are labelled accordingly. Note that each $f(x)\in End(E_{x})$, $f\in C(X)$ are diagonalizable since they commute with their adjoint $f(x)^{*}\in C(X)$. And since such operators commute with each other, the collection $\left\{f(x)\in End(E_{x})|f\in C(X)\right\}$ is simultaneously diagonalizable. Over a point $x\in U_{1}$, the fiber $E_{x}$ decomposes into joint eigenspaces of $\left\{f(x)\in End(E_{x})|f\in C(X)\right\}$. The dimension of these eigenspaces are determined by the number of intersections of the vertical dotted line through $x\in U_{1}$ with the graphs of $\varphi^{1}_{j}$. Using this eigenspace decomposition, we have the following proposition which describes geometrically the right $C(X)$-module structure of $H_{L}$.

\begin{prop}\label{P4.2}
Given $\sigma\in \Gamma(X,E)$ and $f\in C(X)$ the section $\sigma\cdot f\in\Gamma(X,E)$ is given as

\begin{eqnarray}
\left(\sigma\cdot f\right)(x)=\sum\limits_{j=1}^{n} f(\varphi^{i}_{j}(x))e_{j}\cdot\sigma(x)
\label{eq4.1}.
\end{eqnarray}

\noindent where $x\in U_{i}$ and $\sigma(x)=\sum\limits_{j=1}^{n} e_{j}\cdot\sigma(x)$.
\end{prop}

\begin{rem}\label{R4.2}
\begin{enumerate}
\item[]

\item[(1)] In case $C(X)$ is central in $H_{L}$, the above picture reduce to $\left\{U_{i}|i\in I\right\}$ the trivial cover and $\varphi:X\longrightarrow X$ is the identity, i.e. the graph in the above picture is the diagonal of $X\times X$. The action defined by equation (\ref{eq4.1}) then reduces to pointwise multiplication which then coincides with the left $C(X)$-module structure of $H_{L}\cong\Gamma(X,E)$.

\item[(2)] One can understand the right action above as \textit{pointwise-eigenvalue-scaled} action. Compared to the central case, every $f\in C(X)$ acts on a $\sigma\in\Gamma(X,E)$ in a way that $f(x)$ acts diagonally on $\sigma(x)$, i.e. $E_{x}$ constitutes a single eigenspace for the operator $f(x)$ corresponding to the eigenvalue $f(x)\in\mathbb{C}$. In the noncentral case, the action is still pointwise. However, the operator $f(x)$ no longer has a single eigenspace. The eigenspaces are labelled by the points $\varphi^{i}_{j}(x)\in X$ where $x\in U_{i}$ and the eigenvalues of $f(x)$ are $f\left(\varphi^{i}_{j}(x)\right)$, $j=1,...,n$.
\end{enumerate}
\end{rem}

\begin{prop}\label{P4.3}
As a $C(X)$-bimodule, $H_{L}\cong\Gamma(Z,\mathcal{E})$ where $\mathcal{E}$ is a sheaf of complex vector spaces over $X\times X$ supported on a closed subset $Z\subset X\times X$. The $C(X)$-bimodule structure on $\Gamma(Z,\mathcal{E})$ is defined as

\[ (f\cdot \sigma\cdot g)(x,y)=f(x)\sigma(x,y)g(y) \]

\noindent for $f,g\in C(X)$ and $\sigma\in \Gamma(Z,\mathcal{E})$.
\end{prop}

The $C(X)\otimes C(X)^{op}$ is dense in $C(X\times X)$ thus we can extend the $C(X)\otimes C(X)^{op}$-module structure of $H_{L}$ to a $C(X\times X)$-module structure. Consider the annihilator of $H_{L}$,

\[ Ann(H_{L})=\left\{f\in C(X\times X)| f\cdot \sigma = 0, \text{\ for all \ } \sigma\in B_{L}\right\}. \]

\noindent Then, there is an open set $U\subset X\times X$ such that $Ann(H_{L})=C(U)$. Then $Z=(X\times X)-U$, the support of the bimodule $H_{L}$.

\begin{prop}\label{P4.4}
The subset $Z\subseteq X\times X$ is completely determined by the $C(X)$-bimodule structure of $H_{L}$. Moreover, $Z$ is the support of $H_{L}\cong\Gamma(X\times X, \mathcal{E})$.
\end{prop}

By proposition (\ref{P4.1}), $Z$ is the union of transverse subsets of $X\times X$ which is individually are unions of graphs of $\varphi^{i}_{j}$. Let

\[ E_{(x,y)}= \bigoplus\limits_{\varphi^{i}_{j}(x)=y}\left(E_{x}\right)_{\varphi^{i}_{j}(x)} \]

\noindent be the fiber of $\mathcal{E}$ over $(x,y)\in Z$, where $\left(E_{x}\right)_{\varphi^{i}_{j}(x)}$ denotes the eigensubspace of $E_{x}$ over the point $\varphi^{i}_{j}(x)$. This defines a sheaf of vector spaces on $X\times X$ supported on $Z$. A section of $\tau\in\Gamma(X,E)$ defines a section $\hat{\tau}\in \Gamma(Z,\mathcal{E})$ whose value at a point $(x,y)$ is

\[ \hat{\tau}(x,y)=
\begin{cases}
    proj_{ij}\tau(x) ,& \text{if } y=\varphi^{i}_{j}(x) \ \text{for some \ } i,j \\
		& \\
    0,              & \text{otherwise,}
\end{cases} \]

\noindent where $proj_{ij}$ denotes the projection $E_{x}\longrightarrow \left(E_{x}\right)_{\varphi^{i}_{j}(x)}$. Conversely, any section $\tau\in\Gamma(Z,\mathcal{E})$ defines a section $\check{\tau}\in\Gamma(X,E)$ by

\[ \check{\tau}(x)=\sum\limits_{y}\tau(x,y). \]

\noindent Now, given $h\in C(X)\otimes C(X)$ we have

\[ h(x,y)=\sum\limits_{k}f_{k}(x)g_{k}(y) \]

\noindent for some $f_{k},g_{k}\in C(X)$. For any $\sigma\in B\cong\Gamma(X,E)$ we have

\begin{eqnarray*} \left(h\cdot \hat{\sigma}\right)(x,y) &=&\sum\limits_{k}f_{k}(x)\hat{\sigma}(x,y)g_{k}(y) \\
&=& \sum\limits_{k}f_{k}(x)g_{k}(\varphi^{i}_{j}(x))e_{j}\left(\sigma(x)\right) \\
&=& proj_{ij}\left(\sum\limits_{k}f_{k}\cdot\sigma\cdot g_{k}\right)(x,y) \\
\end{eqnarray*}

\noindent which shows that $\phantom{A}^{\wedge}:\Gamma(X,E)\longrightarrow \Gamma(Z,\mathcal{F})$, $\tau\mapsto \hat{\tau}$ is a bimodule map whose inverse is the map $\phantom{A}^{\vee}:\Gamma(Z,\mathcal{F})\longrightarrow \Gamma(X,E)$, $\tau\mapsto\check{\tau}$.

Using proposition (\ref{P4.1}), we can relate the vector bundles the Serre-Swan theorem gives when applied to the left and right $C(X)$-module structure of $H_{L}$ as follows.

\begin{prop}\label{P4.5}
Let $E_{1}\stackrel{p_{1}}{\longrightarrow}X$ and $E_{2}\stackrel{p_{2}}{\longrightarrow}X$ be the vector bundles given by the Serre-Swan theorem applied to the finitely-generated projective left and right $C(X)$-module $H_{L}$, respectively. Then $E_{1}$ and $E_{2}$ are the direct-images of the sheaf $\mathcal{E}$ along $\pi_{1}$ and $\pi_{2}$, respectively.
\end{prop}

First, the direct-image of $\mathcal{E}$ along $\pi_{1}$ is easily seen to be a vector bundle and the space of sections $\Gamma(X,(\pi_{1})_{*}\mathcal{E})$ is easily seen to be isomorphic as left $C(X)$-modules to the left $C(X)$-module $\Gamma(Z,\mathcal{E})$. By proposition (\ref{P4.3}), $\Gamma(Z,\mathcal{E})\cong\Gamma(X,E_{1})$ as left $C(X)$-modules. Thus, by corollary 2.8 of \cite{b007} we see that $E_{1}$ and $(\pi_{1})_{*}\mathcal{E})$ are isomorphic as vector bundles. Similar argument works for $E_{2}$.

Let us say more about the nature of the eigenspaces $E_{(x,y)}$, $x,y\in X$ in relation to the subset $Z$.

\begin{prop}\label{P4.6}
\begin{enumerate}
\item[]
\item[(i)] $E_{x}=\bigoplus\limits_{y\in X} E_{(x,y)}$

\item[(ii)] $dim\left( E_{(x,y)}\right)$ is the number of transverse subsets of $X\times X$ contained in $Z$ passing through $(x,y)$, with multiplicities.

\item[(iii)] $dim\left(\bigoplus\limits_{x\in X}E_{(x,y)}\right)=n$ for any $y\in X$.
\end{enumerate}
\end{prop}

\subsection{The geometry of \textit{C(X)}-ring structures}\label{S4.2}

The previous section describes the geometry of $H_{L}$ using its bimodule structure over $C(X)$. But $H_{L}$ has more structure than just being a bimodule. In particular, it is a $C(X)$-ring via the left source map $C(X)\stackrel{s_{L}}{\longrightarrow}H_{L}$. In this section, we will look at what this additional structure contributes to the geometry of $H_{L}$. We will keep the notations of the previous section.

The $C(X)$-ring structure on $H_{L}\cong\Gamma(X,\mathcal{E})$ via the source map $s_{L}$ consists of a pair of $C(X)$-bimodule maps

\[ \Gamma(Z,\mathcal{E})\tens{C(X)}\Gamma(Z,\mathcal{E}) \stackrel{\mu}{\longrightarrow} \Gamma(Z,\mathcal{E}) \]

\[ C(X)\stackrel{\eta}{\longrightarrow}\Gamma(Z,\mathcal{E}) \]

\noindent satisfying the associativity and unitality conditions. For brevity we will write $\eta=s_{L}$.

\begin{figure}[ht!]
\centering
\includegraphics[width=155mm]{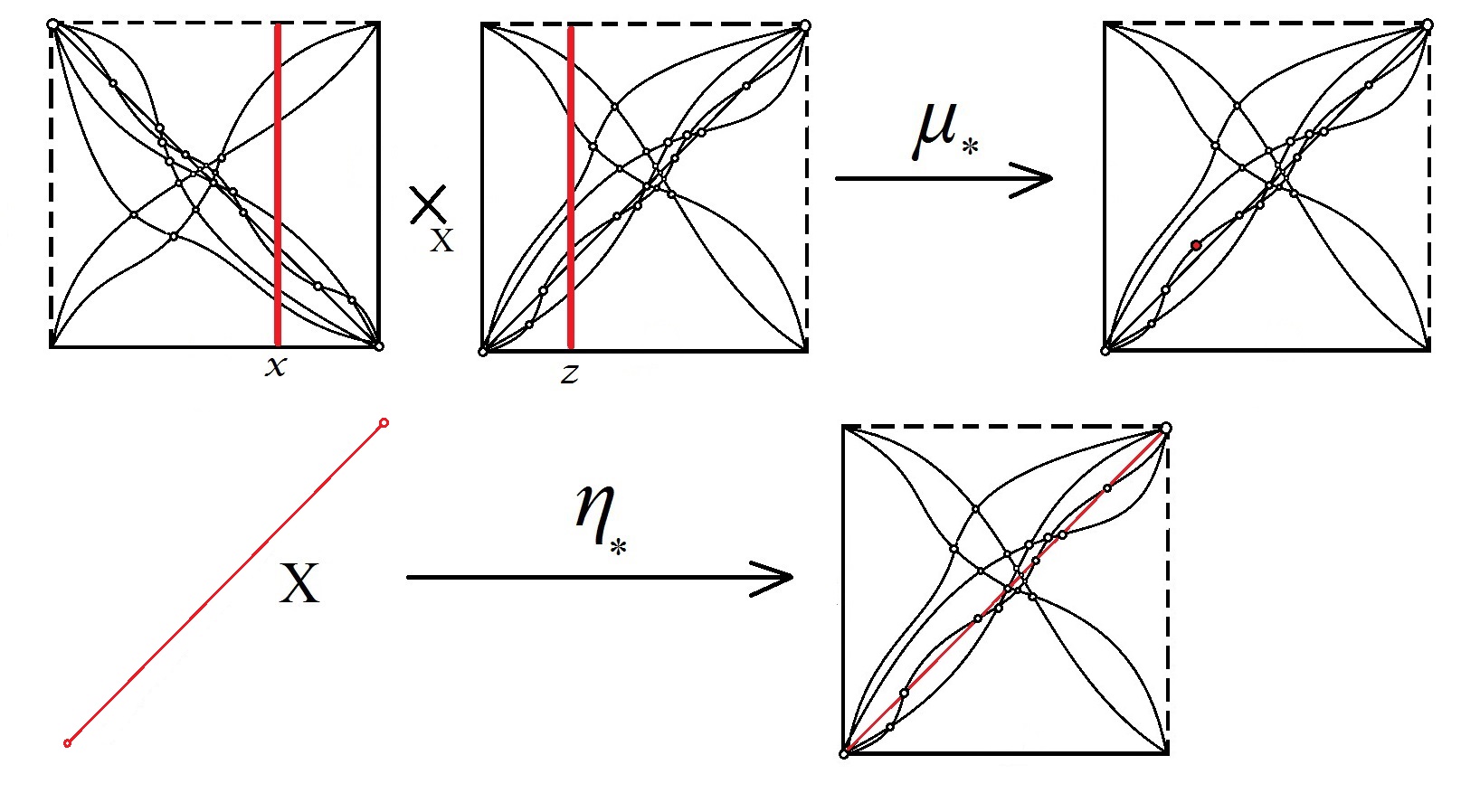}
\caption{The geometry of the product and unit maps.}
\label{fg4.2}
\end{figure}

The unit map $\eta$ gives an element $1\in\Gamma(Z,\mathcal{E})$ satisfying $f\cdot 1=1\cdot f$ for all $f\in C(X)$. Since $X$ is Hausdorff, if $x\neq y$ then we can find an $f\in C(X)$ such that $f(x)=1$ and $f(y)=0$. Thus, for $x\neq y$ we have

\[ 1(x,y)=f(x)1(x,y)=(f\cdot 1)(x,y)=(1\cdot f)(x,y)=1(x,y)f(y)=0. \]

\begin{figure}[ht!]
\centering
\includegraphics[width=135mm]{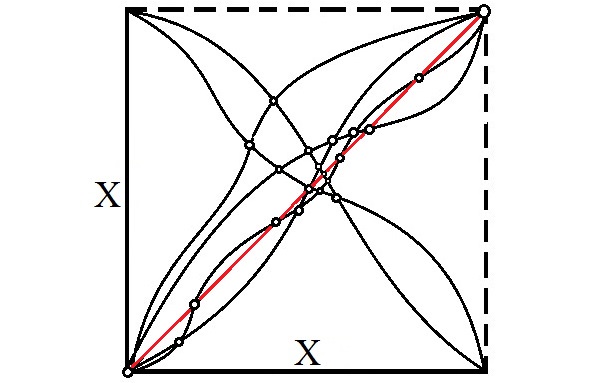}
\caption{Support $Z$ of the bimodule $B$.}
\label{fg4.3}
\end{figure}

\noindent Thus, the source map $A\stackrel{s_{L}}{\longrightarrow} H_{L}$ is implemented by $C(X)\longrightarrow \Gamma(Z,\mathcal{E})$, $f\mapsto f\cdot 1$. This means that $s_{L}\circ f(x)=f(x)1(x,x)$ and choosing $f$ such that $f(x)\neq 0$ and $s_{L}\circ f(x)\neq 0$ we see that $1(x,x)\in E_{(x,x)}$ is a nonzero element. Thus, the diagonal $\slashed{\Delta}$ of $X\times X$ is in $Z$. See figure (\ref{fg4.3}).

Note that $\Gamma(Z,\mathcal{E})\tens{C(X)}\Gamma(Z,\mathcal{E})\cong\Gamma(Z,\mathcal{E}^{(2)})$ where $\mathcal{E}^{(2)}$ is the sheaf of vector spaces whose fiber at a point $(x,z)\in Z$ is the vector space

\[ \bigoplus\limits_{y\in X} \left(E_{(x,y)} \otimes E_{(y,z)}\right) \]

\noindent due to the balancing condition $\sigma\cdot f\otimes_{C(X)}\tau=\sigma\otimes_{C(X)}f\cdot\tau$ for $\sigma,\tau\in\Gamma(Z,\mathcal{E})$ and $f\in C(X)$. Notice that all but finitely many summands above are zero. Specifically, only those $y\in X$ for which $(x,y)$ and $(y,z)$ are both in $Z$ contribute nontrivially. Let us denote these $y\in X$ as $y_{1},y_{2},...,y_{n}$.

By proposition (\ref{P4.3}), $\Gamma(Z,\mathcal{E})\cong\Gamma(X,E)$ as $C(X)$-bimodules. Since $\Gamma(X,-)$ is a fully faithful functor by corollary 2.8 of \cite{b007}, we can convert the global ring structures $\mu$ and $\eta$ into something fiber-wise. In particular, the product map $\mu$ induces a map

\begin{eqnarray}\label{eq4.2}
\xymatrix{ E_{(x,y_{1})}\otimes E_{(y_{1},z)} \oplus ... \oplus E_{(x,y_{n})}\otimes E_{(y_{n},z)} \ar[rr]^-{\mu_{*}} && E_{(x,z)} }
\end{eqnarray}

\noindent illustrated in figure (\ref{fg4.2}). By the universal property of direct sums, there are maps

\[ \xymatrix{ E_{(x,y_{i})}\otimes E_{(y_{i},z)}\ar[rr]^-{\mu^{y_{i}}_{*}} && E_{(x,z)} }\]

\noindent one for each $y_{i}$. The collection of these maps satisfy a set of conditions which, though derivable from associativity, is complicated to write down. See (3) of the remark below for these conditions. However, for the maps $E_{(x,x)}\otimes E_{(x,x)}\stackrel{\mu^{x}_{*}}{\longrightarrow}E_{(x,x)}$ these conditions are precisely the associativity condition. Likewise, the map $\eta$ induces maps $\eta_{*}^{x,y}:\mathbb{C}\longrightarrow E_{(x,y)}$ which is nonzero when $x=y$ and zero otherwise. The map $\mu_{*}^{x}$ together with $\eta_{*}^{x}=\eta_{*}^{x,x}$ makes the vector space $E_{(x,x)}$ a unital algebra, whose dimension depend on the multiplicity of the associated eigenvalue. The following proposition is then immediate from these arguments.

\begin{prop}\label{P4.7}
Let $A^{'}$ be the $C(X)$-sub-bimodule of $\Gamma(Z,\mathcal{E})$ supported on the diagonal $\slashed{\Delta}$. Then $A^{'}$ is an $A$-subring of $H_{L}$ where the multiplication is pointwise. Moreover, $A^{'}$ is the centralizer of $A$ in $H_{L}$.
\end{prop}

\begin{rem}\label{R4.3}
\begin{enumerate}
\item[]

\item[(1)] Using abuse of notation, let us identify $A$ with its image in $H_{L}$. In case $A$ is central in $H_{L}$, the fibers of the vector bundle $E\longrightarrow X$ are algebras. These algebras correspond to $E_{(x,x)}$ together with the maps $E_{(x,x)}\otimes E_{(x,x)}\stackrel{\mu_{*}^{x}}{\longrightarrow} E_{(x,x)}$ and $\mathbb{C}\stackrel{\eta_{*}^{x}}{\longrightarrow} E_{(x,x)}$ since in the central case, $E_{(x,x)}=E_{x}$. Thus, $A^{'}=H_{L}$ in the central case which is not surprising at all knowing that $A^{'}$ is the centralizer of $A$.

\item[(2)] The maps $E_{(x,y_{i})}\otimes E_{(y_{1},z)}\stackrel{\mu^{i}_{*}}{\longrightarrow}E_{(x,z)}$ are only restricted by the associativity of $\mu$. Since $\Gamma(Z,\mathcal{E})\cong\Gamma(X,E)$ and $\Gamma(X,-)$ is known to be a fully faithful functor by corollary 2.8 of \cite{b007}, we have

\[ \xymatrix{ \bigoplus\limits_{y_{i},y_{j}} \left(E_{(x,y_{i})}\otimes E_{(y_{i},y_{j})} \otimes E_{(y_{j},z)}\right) \ar[rr]^-{\left(\bigoplus\limits_{i}\mu^{i}_{*}\right)\otimes id} \ar[dd]_-{id\otimes \left(\bigoplus\limits_{j}\mu^{j}_{*}\right)} && \bigoplus\limits_{y_{j}}\left(E_{(x,y_{j})} \otimes E_{(y_{j},z)}\right) \ar[dd]^-{\bigoplus\limits_{j}\mu^{j}_{*}} \\
 && \\
\bigoplus\limits_{y_{i}}\left(E_{(x,y_{i})} \otimes E_{(y_{i},z)}\right) \ar[rr]_-{\bigoplus\limits_{i}\mu^{i}_{*}} && E_{(x,z).} } \]

\noindent Universal property of direct sums gives us

\[ \xymatrix{ E_{(x,y_{i})}\otimes E_{(y_{i},y_{j})} \otimes E_{(y_{j},z)} \ar[rr]^-{\mu^{i}_{*}\otimes id} \ar[dd]_-{id\otimes \mu^{j}_{*}} && E_{(x,y_{j})} \otimes E_{(y_{j},z)} \ar[dd]^-{\mu^{j}_{*}} \\
 && \\
E_{(x,y_{i})} \otimes E_{(y_{i},z)} \ar[rr]_-{\mu^{i}_{*}} && E_{(x,z)}. } \]

\noindent This justifies the argument before proposition (\ref{P4.7}). We can also use this to say more about the fibers of $\mathcal{E}$ which we state in the next proposition.

\end{enumerate}
\end{rem}

\begin{prop}\label{P4.8}
$E_{(x,y)}$ is a left $E_{(x,x)}-E_{(y,y)}-$bimodule for every $x,y\in X$.
\end{prop}

\begin{rem}\label{R4.4}
Using remark (\ref{R4.3}) above, we can construct a small category $\mathscr{H}_{L}$ enriched over the category of complex vector spaces. The set of objects of $\mathscr{H}_{L}$ is $X$. For every $x,y\in X$, we define

\[  Hom(x,y):=
\begin{cases}
    E_{(x,y)} ,& \text{if } y=\varphi^{i}_{j}(x) \ \text{for some \ } i,j \\
		& \\
    \left\{0\right\},              & \text{otherwise.}
\end{cases} \]

\noindent We will call $\mathscr{H}_{L}$ the \textit{associated category} of the left $A$-bialgebroid $H_{L}$. In the next section, we will see the additional properties of $\mathscr{H}_{L}$ coming from the $A$-coring structure of $H_{L}$. On a different note, let us give a complete geometric description of the $A$-ring structure of $H_{L}$.
\end{rem}

\begin{prop}\label{P4.9}
Denote by $a \ast_{i}b:=\mu_{*}^{y_{i}}(a,b)$, $a\in E_{(x,y_{i})}$ and $b\in E_{(y_{i},z)}$. The product of $\sigma,\tau\in\Gamma(Z,\mathcal{E})$ takes the form

\[ (\sigma\tau)(x,z)=\sum\limits_{i} \sigma(x,y_{i})\ast_{i}\tau(y_{i},z) \]

\noindent for all $(x,z)\in Z$.
\end{prop}

This follows immediately from equation (\ref{eq4.2}). Notice the resemblance of this formula to the one for matrix multiplication. This should remind the reader of an example we discussed in section (\ref{S3.2}). One can view a $C(X)$-ring to be a "matrix" of vector spaces whose entries are indexed by $X\times X$ and what sits in entry $(x,y)$ is the vector space $E_{(x,y)}$. As we have defined after proposition (\ref{P4.4}), the vector space $E_{(x,y)}$ is the zero vector space if $(x,y)\notin Z$. For matrix algebras $M_{n}(\mathbb{C})$, $X$ would be an $n$-element set and the vector spaces $E_{(x,y)}$ would all be $\mathbb{C}$. There are a plethora of algebraic structures package into a bialgebroid let alone in a Hopf algebroid. Before we end this section, let us take a detour to describe the relationships among the structures of $H$: being a $\mathbb{C}$-algebra, the $A$-ring and the $A^{e}$-ring structures being a left-bialgebroid over $A=C(X)$.

For the purpose of this discussion, let us denote by $(\mu_{\mathbb{C}},\eta_{\mathbb{C}})$ the $\mathbb{C}$-algebra structure of $H$ and recall that $(\mu_{L},s_{L})$ and $(\mu_{A^{e}},\eta_{L})$ denote the relevant $A$-ring and $A^{e}$-ring structures of $H$, respectively. As we mentioned in section (\ref{S2.1}), for a $k$-algebra $R$, $R$-ring structures are in bijection with $k$-algebra maps $\eta:k\longrightarrow R$. Thus, the complex algebra structure of $H$ is uniquely determined by the unit map $\eta_{\mathbb{C}}:\mathbb{C}\longrightarrow H$. Similarly, the $A$-ring and the $A^{e}$-ring structures are determined by the $\mathbb{C}$-linear maps $s_{L}$ and $\eta_{L}$. These maps satisfy the following commutativity relations.

\[\xymatrix{\mathbb{C} \ar[rr]^-{\eta_{\mathbb{C}}} \ar[dd] && H\\
&&\\
A \ar[rr] \ar[rruu]|-{s_{L}} && A^{e} \ar[uu]_-{\eta_{L}}} \hspace{.5in}
\xymatrix{H\otimes H \ar@{->>}[rr] \ar[rrdd]_-{\mu_{\mathbb{C}}} && H\tens{A}H \ar@{->>}[rr] \ar[dd]_-{\mu_{L}} && H\tens{A^{e}}H. \ar[lldd]^-{\mu_{A^{e}}} \\
&&&& \\
&& H &&}\]

\noindent In terms of the local eigenspace decomposition, the map $\mu_{\mathbb{C}}$ induces maps

\[\xymatrix{ E_{(x,w)} \otimes E_{(z,y)} \ar[rr] && E_{(x,y)} } \]

\noindent while, by (\ref{eq4.2}), we have maps

\[ \xymatrix{ E_{(x,z)} \otimes E_{(z,y)} \ar[rr] && E_{(x,y)} }. \]

\noindent On the other hand, because the $C(X\times X)$-bimodule structure of $H$ is given as follows,

\[ (f\cdot\sigma)(x,y)=f(x,y)\sigma(x,y), (\sigma\cdot f)(x,y)=f(y,x)\sigma(x,y), \]

\noindent for any $f\in C(X\times X)$, $\sigma\in H$, and $x,y\in X$, the product $\mu_{A^{e}}$ induces maps

\[ \xymatrix{ E_{(x,z)}\otimes E_{(z,x)} \ar[rr] && E_{(x,x)} }. \]

Another way of seeing this is by noting that the product $\mu_{L}$ uses the tensor product $\otimes_{A}$ which kills products $\xymatrix{ E_{(x,w)}\otimes E_{(z,y)} \ar[rr] && E_{(x,y)} }$ for which $w\neq z$. Likewise, the tensor product $\otimes_{A^{e}}$ kills products $\xymatrix{ E_{(x,z)}\otimes E_{(z,y)} \ar[rr] && E_{(x,y)} }$ for which $x\neq y$.

\subsection{The geometry of \textit{C(X)}-coring structures}\label{S4.3}

In this section, using the techniques and results we have developed in sections (\ref{S4.1}) and (\ref{S4.2}) we will describe what the coring structure of $H_{L}$ contributes to the geometry of $\mathcal{E}$. We will keep the notations of the previous two sections.

The $C(X)$-bimodule structure of the underlying $A$-coring structure of $H_{L}$ is related to the $C(X)$-bimodule structure of the underlying $A$-ring via

\begin{equation}\label{eq4.3}
(f\cdot \sigma\cdot g)(x,y)=f(x)g(x)\sigma(x,y)
\end{equation}

\noindent for $\sigma\in \Gamma(Z,\mathcal{E})$, $f,g\in C(X)$, and $x,y\in X$. The left-hand side of equation (\ref{eq4.3}) concerns the bimodule structure one has for the underlying $A$-coring of $H_{L}$ while the right-hand side concerns its $A$-ring structure. This, in particular, implies that if we run the construction we have in section (\ref{S4.1}) for the bimodule structure of the $A$-coring of $H_{L}$, we will get the same sheaf $\mathcal{E}$ supported over the same closed subset $Z$.

The coproduct $\Delta_{L}$ of $H_{L}$, $H_{L}\stackrel{\Delta_{L}}{\longrightarrow}H_{L}\otimes_{A}H_{L}$, uses a different $A$-bimodule structure from the $A$-bimodule structure involved in the $A$-ring structure. Thus, $\otimes_{C(X)}$ means different from the $\otimes_{C(X)}$ we have in the product $\mu$. With this, let us denote by $\boxtimes_{A}$ this new tensor product. thus, we have

\begin{equation}\label{eq4.4}
\xymatrix{\Gamma(Z,\mathcal{E})  \ar[rr]^-{\Delta_{L}} && \Gamma(Z,\mathcal{E})\boxtimes_{C(X)}\Gamma(Z,\mathcal{E}) }.
\end{equation}

\noindent However, using the relation (\ref{eq4.3}) the codomain of $\Delta_{L}$ can be expressed as

\[ \Gamma(Z,\mathcal{E})\boxtimes_{C(X)}\Gamma(Z,\mathcal{E})\cong\Gamma(Z,\mathcal{E}^{\left\langle 2\right\rangle}), \]

\noindent where $\mathcal{E}^{\left\langle 2\right\rangle}$ is the sheaf of vector spaces whose fiber at $(x,z)\in X\times X$ is

\[ \bigoplus\limits_{y^{'},y^{''}\in X}\left(E_{(x,y^{'})}\otimes E_{(x,y^{''})}\right). \]

\noindent Using the same argument we used in the previous section, the map $\Delta_{L}$ induces a map $\left(\Delta_{L}\right)_{*}:\mathcal{E}\longrightarrow\mathcal{E}^{\left\langle 2\right\rangle}$ of sheaves over $Z$. Over point a $(x,y)\in Z$, we have a map

\begin{equation}\label{eq4.5}
\xymatrix{ E_{(x,y)} \ar[rr]^-{\left(\Delta_{L}\right)_{*}^{(x,y)}} && \bigoplus\limits_{z',z''\in X}\left(E_{(x,z')}\otimes E_{(x,z'')}\right) }
\end{equation}

\noindent Meanwhile, the counit $\epsilon_{L}:\Gamma(Z,\mathcal{E})\longrightarrow C(X)$ induces a map $\mathcal{E}\longrightarrow \mathcal{E}^{'}$ of sheaves over $Z$ and $\slashed{\Delta}$, respectively. Here, $Z\longrightarrow\slashed{\Delta}$ is the map $(x,y)\mapsto (x,x)$ for any $(x,y)\in Z$ and $\mathcal{E}^{'}$ is the subsheaf of $\mathcal{E}$ where the fiber of $\mathcal{E}^{'}$ at $(x,y)$ is $\left\{0\right\}$ unless $x=y$, to which the fiber is $\mathbb{C}$ viewed as the one-dimensional subalgebra of $E_{(x,x)}$ spanned by its unit $1(x,x)$. Hence, over a point $(x,y)\in Z$ we have $\left(\epsilon_{L}\right)_{*}^{(x,y)}:E_{(x,y)}\longrightarrow\mathbb{C}$.

Counitality of $\Delta_{L}$ with respect to $\epsilon_{L}$ implies that for fixed but arbitrary $x,y\in X$ we have

\begin{equation}\label{eq4.6} \xymatrix@R=2mm{ & & \bigoplus\limits_{z,z^{'}}\left(E_{(x,z)}\otimes E_{(x,z^{'})}\right) \ar[ddddd]^-{\bigoplus\limits_{z^{'}}id\otimes\left(\epsilon_{L}\right)_{*}^{(x,z^{'})}} \\
& & \\
& & \\
& & \\
& & \\
E_{(x,y)} \ar[rruuuuu]^-{\left(\Delta_{L}\right)_{*}^{(x,y)}} \ar@{=}[rr] & & \bigoplus\limits_{z}\left(E_{(x,z)}\otimes\mathbb{C}\right) \\
v \ar@{|->}[rr] && v\otimes 1 \\} \hspace{.15in}
\xymatrix@R=2mm{ & & \bigoplus\limits_{z,z^{'}}\left(E_{(x,z)}\otimes E_{(x,z^{'})}\right) \ar[ddddd]^-{\bigoplus\limits_{z}\left(\epsilon_{L}\right)_{*}^{(x,z)}\otimes id} \\
& & \\
& & \\
& & \\
& & \\
E_{(x,y)} \ar[rruuuuu]^-{\left(\Delta_{L}\right)_{*}^{(x,y)}} \ar@{=}[rr] & & \bigoplus\limits_{z^{'}}\left(\mathbb{C}\otimes E_{(x,z^{'})}\right) \\
v \ar@{|->}[rr] && 1\otimes v \\}
\end{equation}

\noindent The bottom isomorphisms imply that $\left(\epsilon_{L}\right)_{*}^{(x,z)}$ and $\left(\epsilon_{L}\right)_{*}^{(x,z^{'})}$ are nonzero maps for $z=y$ and $z^{'}=y$. Since $x$ and $y$ are arbitrary to start with, we have the following proposition.

\begin{prop}\label{P4.10}
For any $(x,y)\in Z$, we have $\left(\epsilon_{L}\right)_{*}^{(x,y)}\neq 0$.
\end{prop}

\noindent Another thing we can infer from the diagrams (\ref{eq4.6}), using the isomorphisms in the bottom and the fact that $y$ is among the $z$ and $z^{'}$ that appears as indices, is that the image of $\left(\Delta_{L}\right)_{*}^{(x,y)}$ is contained in

\[  \left(E_{(x,y)}\otimes E_{(x,y)}\right) \oplus \bigoplus\limits_{z,z^{'}} \left(ker\left(id\otimes\left(\epsilon_{L}\right)_{*}^{(x,z^{'})}\right)+ker\left(\left(\epsilon_{L}\right)_{*}^{(x,z)}\otimes id\right)\right) \]

\noindent We will show in the next section that more can be said. In fact, the image of $\left(\Delta_{L}\right)_{*}^{(x,y)}$ is completely contained in $E_{(x,y)}\otimes E_{(x,y)}$.

\subsection{Hopf algebroids over \textit{C(X)}}\label{S4.4}

In this section, we will complete our description of the geometry of the Hopf algebroid $\mathcal{H}$ over $C(X)$. In doing so, we will be able to illustrate the main point of this article. That to such a Hopf algebroid, one can associate a highly structured category.

So far, we have considered only the constituent left bialgebroid $H_{L}$ of $\mathcal{H}$. Running the arguments we have presented in sections (\ref{S4.1}) and (\ref{S4.2}) for $H_{R}$, we see that there is a sheaf of vector spaces $\mathcal{E}^{'}$ over $X\times X$ such that $H_{R}\cong\Gamma(X\times X,\mathcal{E}^{'})$. Let us denote by $Z^{'}$ the support of $H_{R}$ under the isomorphism $H_{R}\cong\Gamma(X\times X,\mathcal{E}^{'})$. The following proposition relates these two sheaves.

\begin{prop}\label{P4.11}
For $H_{L}\cong\Gamma(X\times X,\mathcal{E})$ and $H_{R}\cong\Gamma(X\times X,\mathcal{E}^{'})$ as $C(X)$-bimodules as constructed in sections (\ref{S4.1}) and (\ref{S4.2}), where $\mathcal{E}$ and $\mathcal{E}^{'}$ are sheaves of vector spaces supported on $Z,Z^{'}\subseteq X\times X$, we have
\begin{enumerate}
\item[(1)] $Z=Z^{'}$.
\item[(2)] $\mathcal{E}\cong\mathcal{E}^{'}$ as sheaves over $Z$.
\end{enumerate}
\end{prop}

\begin{prf}
Condition (c) of the definition of a Hopf algebroid implies that the antipode $S$ of $\mathcal{H}$ flips the $C(X)$-bimodule structure used for the $C(X)$-ring structure of $H_{L}$ to that of the $C(X)$-bimodule structure used for the $C(X)$-ring structure of $H_{R}$. Likewise, $S$ flips the bimodule structures of the underlying $C(X)$-coring structures of $H_{L}$ and $H_{R}$. In particular, this tells us that $S$ induces a map $S_{*}:\mathcal{E}\longrightarrow\mathcal{E}^{'}$ which on fibers does $S_{*}(E_{(x,y)})=E_{(y,x)}$ for any $(x,y)\in Z$. Symmetrically, we also have a map denoted the same, $S_{*}:\mathcal{E}^{'}\longrightarrow\mathcal{E}$, which on fibers does $S_{*}(E_{(y,x)})=E_{(x,y)}$ for any $(y,x)\in Z^{'}$. This proves proposition (\ref{P4.11}). $\blacksquare$
\end{prf}

In view of proposition (\ref{P4.11}), we have $\Delta_{R}:\Gamma(Z,\mathcal{E})\longrightarrow\Gamma(Z,\mathcal{E})$. Similar to equation (\ref{eq4.5}), $\Delta_{R}$ induces maps

\begin{equation}\label{eq4.7}
\xymatrix{ E_{(x,y)} \ar[rr]^-{\left(\Delta_{R}\right)_{*}^{(x,y)}} && \bigoplus\limits_{z',z''\in X}\left(E_{(z',y)}\otimes E_{(z'',y)}\right) }
\end{equation}

\noindent for $(x,y)\in Z$. As we promised at the end of section (\ref{S4.3}), $\left(\Delta_{L}\right)_{*}^{(x,y)}$ maps $E_{(x,y)}$ into $E_{(x,y)}\otimes E_{(x,y)}$, for any $(x,y)\in Z$. same holds for $\left(\Delta_{R}\right)_{*}^{(x,y)}$. Let us summarize these statements into the following proposition.

\begin{prop}\label{P4.12}
For every $(x,y)\in Z$,
\begin{enumerate}
\item[(1)] $E_{(x,y)}$ is a coalgebra with coproduct $\left(\Delta_{L}\right)_{*}^{(x,y)}$ and counit $\left(\epsilon_{L}\right)_{*}^{(x,y)}$, and
\item[(2)] $E_{(x,y)}$ is a coalgebra with coproduct $\left(\Delta_{R}\right)_{*}^{(x,y)}$ and counit $\left(\epsilon_{R}\right)_{*}^{(x,y)}$.
\end{enumerate}
\end{prop}

\begin{prf}
We will only prove part $(2)$. The proof for part $(1)$ is similar. The second commutation relation of $\Delta_{L}$ and $\Delta_{R}$ in part $(b)$ of the definition (\ref{D2.2}) gives the following diagram

\begin{equation}\label{eq4.8}
\xymatrix{
E_{(x,y)} \ar[rrrr]^-{\left(\Delta_{L}\right)_{*}^{(x,y)}} \ar[ddd]_-{\left(\Delta_{R}\right)_{*}^{(x,y)}} &&&& \bigoplus\limits_{z^{'},z^{''}}\left(E_{(x,z^{'})}\otimes E_{(x,z^{''})}\right) \ar[dd]^-{\bigoplus\limits_{z^{'}}\left(\Delta_{R}\right)_{*}^{(x,z^{'})}\otimes id}\\
&&&& \\
&&&& \bigoplus\limits_{z^{'},z^{''}}\bigoplus\limits_{\alpha^{'},\alpha^{''}}\left(E_{(\alpha^{'},z^{'})}\otimes E_{(\alpha^{''},z^{'})}\otimes E_{(x,z^{''})}\right) \ar@{=}[d] \\
\bigoplus\limits_{\beta^{'},\beta^{''}}\left(E_{(\beta^{'},y)}\otimes E_{(\beta^{''},y)}\right) \ar[rrrr]_-{\bigoplus\limits_{\beta^{''}}id\otimes\left(\Delta_{L}\right)_{*}^{(\beta^{''},y)}} &&&& \bigoplus\limits_{\beta^{'},\beta^{''}}\bigoplus\limits_{\gamma^{'},\gamma^{''}}\left(E_{(\beta^{'},y)}\otimes E_{(\beta^{''},\gamma^{'})}\otimes E_{(\beta^{''},\gamma^{''})}\right)\\ }
\end{equation}

\noindent for fixed but arbitrary $(x,y)\in Z$. In the composite $\left(\bigoplus\limits_{z^{'}}\left(\Delta_{R}\right)_{*}^{(x,z^{'})}\otimes id\right)\circ\left(\Delta_{L}\right)_{*}^{(x,y)}$, the third leg lands in $\bigoplus\limits_{z^{''}}E_{(x,z^{''})}$. On the other hand, the third leg of the composite $\left(\bigoplus\limits_{\beta^{''}}id\otimes\left(\Delta_{L}\right)_{*}^{(\beta^{''},y)}\right)\circ\left(\Delta_{R}\right)_{*}^{(x,y)}$ lands in $\bigoplus\limits_{\beta^{''},\gamma^{''}}E_{(\beta^{''},\gamma^{''})}$. This implies that for $\beta^{''}\neq x$, $E_{(\beta^{''},y)}\subseteq ker\ \left(\Delta_{L}\right)_{*}^{(\beta^{''},y)}$. From our last statement in section (\ref{S4.3}), $\left(\Delta_{L}\right)_{*}^{(\beta^{''},y)}\left(E_{(\beta^{''},y)}\right)$ is contained in

\[  \left(E_{(\beta^{''},y)}\otimes E_{(\beta^{''},y)} \right)\oplus\bigoplus\limits_{f^{'},f^{''}}\left(ker\left(id\otimes\left(\epsilon_{L}\right)_{*}^{(\beta^{''},f^{'})}\right)+ker\left(\left(\epsilon_{L}\right)_{*}^{(\beta^{''},f^{''})}\otimes id\right)\right). \]

\noindent Counitality of $\Delta_{L}$ with respect to $\epsilon_{L}$, implemented locally by diagram (\ref{eq4.6}), gives

\[ \xymatrix@R=.2cm{ E_{(\beta^{''},y)} \ar[r]^-{\cong}  & \bigoplus\limits_{f^{'}}\left(id\otimes\left(\epsilon_{L}\right)_{*}^{(\beta^{''},f^{'})}\right)\left(\Delta_{L}\right)_{*}^{(\beta^{''},y)}(E_{(\beta^{''},y)}) \ar@{=}[r] & \left\{0\right\}. \\
v \ar@{|->}[r] & v\otimes 1 & \\} \]

\noindent By assumption, $E_{(\beta^{''},y)}$ are nontrivial. This is a contradiction unless the summands corresponding to $\beta^{''}\neq x$ of the direct sum in the lower left corner of diagram (\ref{eq4.8}) do not intersect the image of $\left(\Delta_{R}\right)_{*}^{(x,y)}$.

Using the first commutation relation in part $(b)$ of the definition (\ref{D2.2}), we get a diagram similar to diagram (\ref{eq4.8}). Inspecting that resulting diagram tells us that the image of $\left(\Delta_{R}\right)_{*}^{(x,y)}$ does not intersect those summands of the direct sum in the lower left corner of diagram (\ref{eq4.8}) corresponding to $\beta^{'}\neq x$. This shows that, indeed,

\[ \left(\Delta_{R}\right)_{*}^{(x,y)}: E_{(x,y)}\longrightarrow E_{(x,y)}\otimes E_{(x,y)}. \]

\noindent The coassociativity of $\left(\Delta_{R}\right)_{*}^{(x,y)}$ follows from coassociativity of $\Delta_{R}$ and its counitality with respect to $\left(\epsilon_{R}\right)_{*}^{(x,y)}$ follows from counitality of $\Delta_{R}$ with respect to $\epsilon_{R}$. This proves part $(2)$ of the above proposition. Exchanging the roles of $\Delta_{L}$ and $\Delta_{R}$ with minor modifications proves part $(1)$. $\blacksquare$
\end{prf}

Following the arguments in sections (\ref{S4.1}), (\ref{S4.2}) and (\ref{S4.3}) for $H_{R}$, we see that we can similarly associate a category $\mathscr{H}_{R}$ enriched over $\mathcal{V}$. Denoting by $C(\mathcal{V})$ by the category of coalgebras on $\mathcal{V}$, we have the following proposition.

\begin{prop}\label{P4.13}
The categories $\mathscr{H}_{L}$ and $\mathscr{H}_{R}$ are enriched over $C(\mathcal{V})$.
\end{prop}

These categories are strongly related. By proposition (\ref{P4.11}), we have the following corollary.

\begin{cor}\label{P4.14}
The $C(\mathcal{V})$-enriched categories $\mathscr{H}_{L}$ and $\mathscr{H}_{R}$ have isomorphic underlying $\mathcal{V}$-enriched categories.
\end{cor}

Note that the underlying $\mathcal{V}$-enriched category of $\mathscr{H}_{L}$ and $\mathscr{H}_{R}$ only depends on the $C(X)$-ring structures of $H_{L}$ and $H_{R}$, respectively. Another way to prove corollary (\ref{P4.14}) is to use the fact that $H_{L}$ and $H_{R}$ have the isomorphic $C(X)$-ring structures. To see why $H_{L}$ and $H_{R}$ have isomorphic $C(X)$-ring structures, note that the source map of $H_{L}$ is the target map of $H_{L}$ while the target map of $H_{L}$ is the source map of $H_{R}$. In the general definition of a Hopf algebroid, one can either use the source or the target map to select a particular ring structure to consider, see for example \cite{bohm}. Using the general fact that for a general $k$-algebra $R$, $R$-rings $(A,\mu,\eta)$ corresponds uniquely to $k$-algebra maps $\eta$, we see that $H_{L}$ and $H_{R}$ are isomorphic as $C(X)$-rings.

\begin{rem}\label{R4.5}
Another way to see why $H_{L}$ and $H_{R}$ are isomorphic as $C(X)$-rings is the fact that general Hopf algebroids $\mathcal{H}$ with bijective antipode over a \textit{commutative} ring $K$ is a coupled $K$-Hopf algebra.
\end{rem}

Unlike the ring structures, the $C(X)$-coring structures of $H_{L}$ and $H_{R}$ can vary wildly as illustrated by coupled Hopf algebras. This implies that the $C(\mathcal{V})$-enrichments $\mathscr{H}_{L}$ and $\mathscr{H}_{R}$ need not be isomorphic. However, they form a topological coupled Hopf category. The coupling functor is the one induced by the antipode $S$ of the Hopf algebroid $\mathcal{H}$. We formalize this in the following theorem.

\begin{thm} \label{T4.1}
Given a finitely-generated projective Hopf algebroid $\mathcal{H}$ over $C(X)$ with bijective antipode, one can associate a topological coupled Hopf category $\mathscr{H}$ via the construction we presented in sections (\ref{S4.1}) and (\ref{S4.2}). Conversely, to any topological coupled Hopf category $\mathscr{H}$, the space of sections $\Gamma(X\times X,H)$ of the associated sheaf $H$ of $\mathscr{H}$ is a Hopf algebroid over $C(X)$.
\end{thm}

The proof of the first statement is basically the breadth of section (\ref{S4.0}). For the second statement, one can consider the bimodule structures presented in proposition (\ref{P4.3}). The rest of the structures are given by the rest of the structure maps of $\mathscr{H}$. The above theorem is a generalization of the example in \cite{bcv} where they constructed out of a $k$-linear category with finitely many objects a weak Hopf algebra. The theorem not only recovers an inverse to the construction they presented but it also work for weak Hopf algebra as long as the subalgebra spanned by the left and the right units are commutative. The above theorem is the generalization of the example we discussed in section (\ref{S3.2}).

\subsection{The central case}\label{S4.5}

Although $C(X)$ is commutative, it may not be central in $\mathcal{H}$. Let us look at the special case when $C(X)$ is central in $\mathcal{H}$, by which we mean that the images of the source and target maps are central in the relevant $C(X)$-ring structure of $\mathcal{H}$. For simplicity, we will blur the disctinction between $C(X)$ at its images under the source maps of $\mathcal{H}$.

Let us consider first the constituent left bialgebroid $H_{L}$ of $\mathcal{H}$. By proposition (\ref{P4.7}), $H_{L}$ is supported along the diagonal $\slashed{\Delta}\subseteq X\times X$. This means that the sheaf $\mathcal{E}$ coincides with the vector bundle $E\longrightarrow X$. We can simply identify the diagonal $\slashed{\Delta}$ with $X$. With this, the multiplication $\mu_{L}$ in $H_{L}$ via the identification $H_{L}\cong\Gamma(X,E)$ is pointwise, i.e. the fibers of the vector bundle $E\longrightarrow X$ are (possibly nonisomorphic) unital complex alegrbas $\left(E_{x},\left(\mu_{L}\right)_{*}^{x},\left(s_{L}\right)_{*}^{x}\right)$, where $\left(\mu_{L}\right)_{*}^{x}$ and $\left(s_{L}\right)_{*}^{x}$ are the maps induced by $\mu_{L}$ and $s_{L}$ on the fiber $E_{x}$.

By proposition (\ref{P4.12}), the coproduct $\Delta_{L}$ and counit $\epsilon_{L}$ of $H_{L}$ also descends into a coproduct $\left(\Delta_{L}\right)_{*}^{x}$ and a counit $\left(\epsilon_{L}\right)_{*}^{x}$ for the fibers $E_{x}$, $x\in X$, making them coalgebras. Using condition $(b)$ in the definition of a bialgebroid, we see that $\left(\Delta_{L}\right)_{*}^{x}$ is multiplicative for any $x\in X$. Meanwhile, using condition $(c)$ of the definition of a bialgebroid we see that $\left(\epsilon_{L}\right)_{*}^{x}$ is multiplicative for any $x\in X$. This gives us the following proposition.

\begin{prop}\label{P4.15}
If $C(X)$ is central in $H_{L}$, then for any $x\in X$,

\[ \left(E_{x},\left(\mu_{L}\right)_{*}^{x},\left(s_{L}\right)_{*}^{x},\left(\Delta_{L}\right)_{*}^{x},\left(\epsilon_{L}\right)_{*}^{x}\right) \]

\noindent is a bialgebra. Moreover, the bialgebroid $H_{L}$ is a bundle of bialgebras via $H_{L}\cong\Gamma(X,E)$.
\end{prop}

Similar statement holds for the constituent right bialgebroid $H_{R}$. Since for very $x\in X$ the maps $(s_{L})_{*}^{x}$ and $(s_{R})_{*}^{x}$ induced by the source maps $s_{L}$ and $s_{R}$ are the same, the multiplications $(\mu_{L})_{*}^{x}$ and $(\mu_{R})_{*}^{x}$ coincide. Assuming mild nondegeneracy conditions for $(\Delta_{L})_{*}^{x}$ and $(\Delta_{R})_{*}^{x}$, we get the following proposition.

\begin{prop}\label{P4.16}
Let $\mathcal{H}=(H_{L},H_{R},S)$ be a Hopf algebroid over $A=C(X)$ where $A$ is central in both $H_{L}$ and $H_{R}$. Denote by $H$ the underlying complex algebra of $\mathcal{H}$. Suppose that the maps

\[ \xymatrix@R=2mm{H\tens{A}H \ar[rr]^-{\mathfrak{gal}_{L}} && H\tens{A}H \\
a\tens{A}b \ar@{|->}[rr] &&  ab_{[1]}\tens{A}b_{[2]}\\}, \hspace{.5in}
\xymatrix@R=2mm{H\tens{A}H \ar[rr]^-{\mathfrak{gal}_{R}} && H\tens{A}H \\
a\tens{A}b \ar@{|->}[rr] &&  ab^{[1]}\tens{A}b^{[2]}\\} \]

\noindent are bijections. Then
\begin{enumerate}
\item[(i)] $H$ is a coupled Hopf algebra with constituent Hopf algebras $H_{L}$ and $H_{R}$ and coupling map $S$.
\item[(ii)] Each fiber $E_{x}$ is a Hopf algebra and $H_{L}\cong\Gamma(X,E)$ as Hopf algebras, where the structure maps of $\Gamma(X,E)$ are all pointwise. Same is true for $H_{R}$.
\item[(iii)] $\mathcal{H}$ is a bundle of coupled Hopf algebras over $X$ such that the constituent Hopf algebras at a point $x\in X$ are the fiber Hopf algebras of $H_{L}$ and $H_{R}$.
\end{enumerate}
\end{prop}

\begin{prf}
Centrality of $A$ in both $H_{L}$ and $H_{R}$ implies that $H_{L}$ and $H_{R}$ are in fact bialgebras over $A$ (not just bialgebroids). The nondegeneracy conditions assumed in the proposition implies that $H$ is a Galois extension for both bialgebras $H_{L}$ and $H_{R}$. By \cite{sch1}, the bialgebras $H_{L}$ are $H_{R}$ are in fact Hopf algebras, i.e. the identity maps $H_{L}\stackrel{id}{\longrightarrow}H_{L}$ and $H_{R}\stackrel{id}{\longrightarrow}H_{R}$ are invertible in the respective convolution algebras associated to the bialgebras $H_{L}$ and $H_{R}$. The rest of the conditions for $\mathcal{H}$ to be a Hopf algebroid imply that $H_{L}$ and $H_{R}$ are coupled Hopf algebras with coupling map $S$, the antipode of $\mathcal{H}$. This proves part $(i)$.

To prove part $(ii)$, we argue that the maps $\mathfrak{gal}_{L}$ and $\mathfrak{gal}_{R}$ are $A$-bimodule maps. Thus, there descend into fiberwise bijections. Using the same argument we did for part $(i)$, see that the fibers are coupled Hopf algebras. Part $(iii)$ readily follows from the proofs of parts $(i)$ and $(ii)$. $\blacksquare$
\end{prf}

\section{Correspondence of Galois extensions}\label{S5.0}

In this section, we will see that the correspondence between Hopf algebroids and coupled Hopf categories we established in theorem (\ref{T4.1}) persists to their corresponding Galois theories. To be precise, we will prove the following theorem.

\begin{thm} \label{T5.1}
Let $\mathcal{H}=(H_{L},H_{R},S)$ be a Hopf algebroid over $A=C(X)$ for some compact Hausdorff space $X$. Let $\mathscr{H}$ be the corresponding topological coupled Hopf category of $\mathcal{H}$. Then $\mathcal{H}$-Galois extensions of $A$ corresponds bijectively to $\mathscr{H}$-Galois extensions of $\mathbb{1}^{X}$.
\end{thm}

Before proving the above theorem, let us comment on what we mean by Galois extension by a (topological) coupled Hopf category $\mathscr{H}=(\mathscr{H}_{L},\mathscr{H}_{R},S)$. By this, we mean an inclusion of categories $\mathbb{1}^{X}\subseteq \mathscr{M}$ which is simultaneously $\mathscr{H}_{L}$-Galois and $\mathscr{H}_{R}$-Galois in the sense of section (\ref{S3.3}). Note that by definition (\ref{D3.3}), we are not requiring $\mathscr{H}_{L}$ and $\mathscr{H}_{R}$ to be Hopf categories (individually, they are only $C(\mathcal{V})$-enriched categories). In particular, they do not necessarily have antipodes. Fortunately, Galois extension in the sense described in section (\ref{S3.3}) does not really make use of the antipode.

\begin{prf}
Let $B$ be a (left) $\mathcal{H}$-Galois extension of $A$. In particular, $B$ is an $A$-ring. Note that the arguments we used in sections (\ref{S4.1}) and (\ref{S4.2}) only use the $A$-ring structure of the Hopf algebroid $\mathcal{H}$. Using the same arguments, $B\cong\Gamma(X\times X, \mathcal{B})$ where $\mathcal{B}$ is a sheaf of vector spaces over $X\times X$. By the Galois condition, we see that $\mathcal{B}$ has the same support $Z\subseteq X\times X$ as the sheaf $\mathcal{E}$ we get from either $H_{L}$ or $H_{R}$. Similar to remark (\ref{R4.4}), we get a small category $\mathscr{B}$ over $X$ enriched over $\mathcal{V}$ whose associated sheaf is $\mathcal{B}$.

The (right) $H_{L}$-coaction $\rho_{L}:B\longrightarrow B\otimes_{A}H$ induces a map $\mathcal{B}\longrightarrow \mathcal{B}\prescript{}{X}{\times}_{X}\mathcal{E}$ of sheaves of $\mathcal{O}_{X}$-bimodules over $X\times X$. By definition, $B$ is a right $A$-module and a right $A^{op}$-module. Using this, the $A$-bimodule structure on $B$ is as follows:

\[ a\cdot b\cdot a^{'}=b(aa^{'}) \]

\noindent for any $a,a^{'}\in A$ and $b\in B$. Similar to (\ref{eq4.5}), the right $H_{L}$-coaction induces, for every $(x,y)\in Z$, linear maps

\begin{equation}\label{eq5.1}
\xymatrix{ B_{(x,y)} \ar[rr]^-{(\rho_{L})_{*}^{(x,y)}} && \bigoplus\limits_{z^{'},z^{''}\in X} B_{(z^{'},y)}\otimes E_{(z^{''},y)}}
\end{equation}

\noindent where $B_{(x,y)}$ the fiber of $\mathcal{B}$ at the point $(x,y)$. As before, $E_{(x,y)}$ denotes the fiber of $\mathcal{E}$ over $(x,y)$. Likewise, the right $H_{R}$-coaction $\rho_{R}$ induces linear maps

\begin{equation}\label{eq5.2}
\xymatrix{ B_{(x,y)} \ar[rr]^-{(\rho_{R})_{*}^{(x,y)}} && \bigoplus\limits_{z^{'},z^{''}\in X} B_{(x,z^{'})}\otimes E_{(x,z^{''})}}
\end{equation}

By diagram (\ref{eq2.1}), we have

\begin{equation}\label{eq5.3}
\xymatrix{
B_{(x,y)} \ar[rrrr]^-{\left(\rho_{L}\right)_{*}^{(x,y)}} \ar[ddd]_-{\left(\rho_{R}\right)_{*}^{(x,y)}} &&&& \bigoplus\limits_{z^{'},z^{''}}\left(B_{(x,z^{'})}\otimes E_{(x,z^{''})}\right) \ar[dd]^-{\bigoplus\limits_{z^{'}}\left(\rho_{R}\right)_{*}^{(x,z^{'})}\otimes id}\\
&&&& \\
&&&& \bigoplus\limits_{z^{'},z^{''}}\bigoplus\limits_{\alpha^{'},\alpha^{''}}\left(B_{(\alpha^{'},z^{'})}\otimes E_{(\alpha^{''},z^{'})}\otimes E_{(x,z^{''})}\right) \ar@{=}[d] \\
\bigoplus\limits_{\beta^{'},\beta^{''}}\left(B_{(\beta^{'},y)}\otimes E_{(\beta^{''},y)}\right) \ar[rrrr]_-{\bigoplus\limits_{\beta^{''}}id\otimes\left(\Delta_{L}\right)_{*}^{(\beta^{''},y)}} &&&& \bigoplus\limits_{\beta^{'},\beta^{''}}\bigoplus\limits_{\gamma^{'},\gamma^{''}}\left(B_{(\beta^{'},y)}\otimes E_{(\beta^{''},\gamma^{'})}\otimes E_{(\beta^{''},\gamma^{''})}\right)\\ }
\end{equation}

Meanwhile, counitality of the left coaction $\rho_{L}$ implies that

\[ \xymatrix@R=2mm{ & & \bigoplus\limits_{z^{'},z^{''}}\left(B_{(x,z^{'})}\otimes E_{(x,z^{''})}\right) \ar[ddddd]^-{\bigoplus\limits_{z^{''}}id\otimes\left(\epsilon_{L}\right)_{*}^{(z,y)}} \\
& & \\
& & \\
& & \\
& & \\
B_{(x,y)} \ar[rruuuuu]^-{\left(\rho_{L}\right)_{*}^{(x,y)}} \ar@{=}[rr] & & \bigoplus\limits_{z^{'}}\left(B_{(x,z^{'})}\otimes\mathbb{C}\right) \\
v \ar@{|->}[rr] && v\otimes 1 \\} \]

\noindent from which, using a similar argument we to the proof of proposition (\ref{P4.12})(1), gives

\[ \xymatrix{ B_{(x,y)} \ar[rr]^-{(\rho_{L})_{*}^{(x,y)}} && B_{(x,y)}\otimes E_{(x,y)}}. \]

\noindent Similarly, we have

\[ \xymatrix{ B_{(x,y)} \ar[rr]^-{(\rho_{R})_{*}^{(x,y)}} && B_{(x,y)}\otimes E_{(x,y)}}. \]

\noindent These tell us that $\mathscr{B}$ is a right $\mathscr{H}_{L}$- and a right $\mathscr{H}_{R}$-comodule. The composition $\circ$ in $\mathscr{B}$ is induced by the $A$-product on $B$. By equations (\ref{eq2.2}) to (\ref{eq2.5}), this composition $\circ$ is a map of right $\mathscr{H}_{L}$- and a right $\mathscr{H}_{R}$-modules. Thus, $\mathscr{B}$ is a right $\mathscr{H}_{L}$- and a right $\mathscr{H}_{R}$-comodule-category. It is not hard to see that the right coactions of $\mathscr{H}_{L}$ and $\mathscr{H}_{R}$ on $\mathscr{B}$ are both Galois whose subcategories of coinvariants are both the same as $I_{X}$. These imply that $\mathscr{B}$ is a Galois extension of $I_{X}$ by the topological coupled Hopf category $\mathscr{H}=(\mathscr{H}_{L},\mathscr{H}_{R},S)$.

The inverse of this correspondence is easily seen as the the one that associates to an $(\mathscr{H}_{L},\mathscr{H}_{R},S)$-Galois extension $I_{X}\subseteq \mathscr{B}$ the $(H_{L},H_{R},S)$-Galois extension $A\subseteq B$ where $H_{L},H_{R},B$ and $A$ are the space of global sections of the associated sheaves to $\mathscr{H}_{L},\mathscr{H}_{R},\mathscr{B}$ and $I_{X}$, respectively. The compatibility conditions in the categorical side precisely correspond to the analogous compatibility conditions in the algebraic side. $\blacksquare$
\end{prf}


\end{document}